\documentclass[12pt]{article}
\usepackage{amssymb}
\usepackage{latexsym}
\usepackage{amscd}

\usepackage{amsmath}
\usepackage{a4}
\usepackage{epsf}
\usepackage{graphics}
\def\a{\alpha}
\def\b{\beta}
\def\g{\gamma}
\def\G{\Gamma}
\def\d{\delta}

\def\ape{\varepsilon}
\def\la{\lambda}

\def\pat{\partial}
\def\o{\omega}

\def\r{\rho}
\def\phi{\varphi}
\def\ta{\theta}
\newcommand{\osc}{\operatorname{osc}}

\newcommand{\qed}{\hfill\rule{5pt}{5pt}}
\begin{document}
\title{Conley Index Theory and Novikov-Morse Theory }

\author{Huijun Fan\thanks{The first author is partially supported by
Research Fund for Returned Overseas Chinese\ Scholars 20010107 and
the Partner Group of the Max-Planck-Institute for Mathematics in
the Sciences in Leipzig at the Chinese Academy of Sciences.} \quad
J\"urgen Jost}
\date{}

\maketitle
\begin{abstract}
\noindent We derive general Novikov-Morse type inequalities in a
Conley type framework for flows carrying cocycles, therefore
generalizing our results in \cite{FJ2} derived for integral
cocycle. The condition of carrying a cocycle expresses the
nontriviality of integrals of that cocycle on flow lines.
Gradient-like flows are distinguished from general flows carrying
a cocycle by boundedness conditions on these integrals.
\end{abstract}
{\bf MSC2000:}\quad 37B30,37B35,57R70\\

\vskip 20pt \noindent
{\bf {1. Introduction}}\\
  \\
\indent In the 80's, S. P. Novikov \cite{No1, No2} considered
Morse type inequalities relating the zeros of a closed Morse
1-form $\o$ to the topology of the underlying space $X$. By the
Poincar\'e lemma, a closed 1-form is locally exact. For instance,
we can assume that $\o=df_p$ near a zero point $p$ of $\o$. So
closed Morse 1-form means that each zero point $p$ of $\o$ is a
non-degenerate critical point of $f_p$ and has a Morse index from
$f_p$. Define $S(\o)$ to be the union of the sets $S_i(\o)$ of the
zero points of $\o$ with index $i$. Define
\begin{align*}
c_i:=\# S_i(\o)\\
\end{align*}
Then the Morse type inequalities are
\begin{align*}
c_i\ge & b_i([\o])\\
\sum^i_{j=0}(-1)^{i-j}c_j\ge&
\sum^i_{j=0}(-1)^{i-j}b_j([\o])\tag{$1.1$}
\end{align*}
for $i=0,1\cdots,m$. Here the Novikov numbers $b_i([\o])$ are
determined only by $X$ and the cohomology class $[\o]$. These
Morse type inequalities involving Novikov numbers are called
Novikov inequalities.

Novikov inequalities were obtained by constructing the Novikov
complex $(C_*,\pat)$ with respect to the closed Morse 1-form $\o$.
If $[\o]$ is an integral cohomology class, then each $C_i$ is a
free ${\Bbb Z}((t))$-module generated by the set of all the zeros
of $\o$ with index $i$. Here ${\Bbb Z}((t)):={\Bbb
Z}[[t]][t^{-1}]$ is the Novikov ring. For the detailed
construction of the boundary operators and some of their
properties, see the discussion in \cite{BF}, \cite{Pa}, \cite{Ra}
and the references there.

M. Farber and A.A. Ranicki \cite{FR} have used a noncommutative
localization method to construct a universal complex for a rank
$1$ Morse closed 1-form. Subsequently, M. Farber \cite{Fa2}
applied the method to the general case, which can induce many
kinds of Novikov complexes. In a series of papers \cite{Fo1, Fo2,
Fo3}, R. Forman developed the discrete Morse theory and Novikov
theory in a combinatorial category. Those Novikov-Morse
inequalities are proved also by constructing certain combinatory
Morse(Novikov) complexes, e.g., using the Witten deformation
technique. In order to use the techniques from the smooth
category, he introduced combinatorial vector fields and
differential forms.

Novikov theory for closed 1-forms is very important, since many
functionals, such as the symplectic action functional, the
Chern-Simons functional and many Hamiltonian systems from
electromagnetism and fluid mechanics are multi-valued functionals.
Though these are multi-valued functionals in infinite dimensional
spaces, people hope that the Novikov theory in finite dimension
can provide a model for the infinite dimensional case, and that it
can provide a strong topological method to find the critical
points of the multi-valued functionals.

Though Novikov theory for closed 1-forms has been studied for many
years, the various proofs of Novikov inequalities were all based
on constructing various (analytic or topological )Novikov
complexes. Hence the closed 1-form needs some non-degeneracy
assumptions. One could not handle the case if the closed 1-form
$\o$ has general degenerate zero points ( zero locus). On the
other hand, the behavior of the dynamical systems generated by a
closed 1-form has not been carefully studied. For instance, the
previous research was limited to the smooth category. Though there
should be an analogue in the continuous category, if so, then
comes the natural question, what is the relation between this
analogue and the well-known gradient-like dynamical systems?

In \cite{FJ2}, we introduced a concept of ``flow carrying a
cocycle $\a$ ($\a$-flow and generalized $\a$-flow)'', where $\a$
is a 1-dimensional cocycle in the bounded Alexander-Spanier
cohomology (see i.e., \cite{Sp}). This is the analogue of the flow
generated by a closed 1-form. Actually, as shown in \cite{FJ2},
(generalized) $\a$-flows characterize most important dynamical
systems, in particular the following:
\begin{itemize}
\item The gradient flow of a Morse function $f$ is a $df$-flow.

\item The gradient flow of a closed 1-form $\o$ is an $\o$-flow
and certain perturbation flows are also $\o$-flows.

\item If an $\a$-flow has no fixed point, then this flow is a flow
``carrying a cohomology class'' as introduced in \cite{Ch}.

\item An $\a$-Morse-Smale flow as defined in \cite{FJ2} which is a
generalization of the famous Morse-Smale flow; the generalization
here essentially consists in allowing the existence of "cycles"

\item If $\a$ is a coboundary, then this flow is a gradient-like
flow. Conversely, if a flow is a gradient-like flow ,then this
flow carries a coboundary $\d g$ where $g$ is a Lyapunov function.

\end{itemize}
\indent A flow carrying a cocycle $\a$ has a $\pi$-Morse
decomposition (Theorem 3.4.4 in \cite{FJ2}). If $\a$ is an
integral cocycle, we can obtain the Novikov-Morse type
inequalities, Theorem 4.4.1 in \cite{FJ2}. This theorem uses the
Conley index to
 get the ``local'' topological information of the isolated invariant sets.\\
\indent Those inequalities in \cite{FJ2} are a generalization of
the Novikov inequalities for an integral closed Morse 1-form.
Starting from these inequalities, we can recover many Novikov type
inequalities found before. For example, we can recover the Novikov
inequalities for an integral closed 1-form $\o$ having Bott type
nondegenerate zero sets, which was given in \cite{Fa2}. In
addition, our theory has some new features:
\begin{itemize}

\item Novikov inequalities now hold in the continuous category
(i.e.for (generalized) $\a$-flows).

\item  New Novikov inequalities (Theorem 4.4.1 in \cite{FJ2}) for
a general closed 1-form without any non-degeneracy requirement.

\item  New Novikov inequalities for $\a$-Morse
  Smale flows.

\item  Vanishing of the Novikov numbers $b_i([\a])$, if the
manifold allows the existence of a flow carrying a cohomology class $\a$.\\

\end{itemize}

However, \cite{FJ2} is mainly focused on the explanation of the
flow carrying a cocycle $\a$, and on the proof of the
Novikov-Morse inequalities for this flow under the assumption that
$\a$ is an integral cocycle.

This paper is the second one. There are two aims. The first aim is
to give a sufficient and necessary description of the relation
between $\a$-flows and the gradient-like flows. This will be given
in section 3. The second aim is to generalize the Novikov-Morse
inequalities for integral cocycles $\a$ to the analogous results
for higher rank cocycles $\a$. This is discussed in sections 4,5.
To do this, one should prove an analogue of Theorem 4.2.1 in
\cite{FJ2}. This is a nontrivial generalization, since the exit
set of the lifting flow becomes more complicated when compared
with the integral cocycle case. One should handle carefully the
deformation complex in the corners. As long as we have the
analogous theorem, we can easily generalize those theorems
in \cite{FJ2} by following proofs in \cite{FJ2}.\\
\indent In section 2, we recollect some basic definitions and
examples of flows carrying a cocycle $\a$.

Recently, M. Farber \cite{Fa3}(2001) also proposed the
concept:"Lyapunov one-form" $\o$ for a pair $(\Phi, Y)$, where
$\Phi$ is a continuous flow and $Y$ is a closed invariant set in
$\Phi$. Since a flow $\Phi$ allowing the existence of a Lyapunov
one-form $\o$ is actually a flow carrying a cocycle $\o$
introduced in \cite{FJ2}, our definition seems more general than
his. In \cite{Fa3,Fa4}, Farber constructed the corresponding
Ljusternik-Schnirelman theory. In \cite{FKLZ}, the authors
discussed the existence of "Lyapunov one-forms".

The results presented in this paper were obtained in 2000 and
circulated in preliminary preprint form in \cite{FJ1}.\\
 \\
{\bf  2. Review of flows carrying a cocycle}\\
 \\
\indent In \cite{FJ2}, we have defined the dynamical systems,
flows carrying a cocycle. In this section we will recollect some
facts
about such dynamical systems.\\
\indent Let $(X,d)$ be a compact metric space with metric $d$. Let
$v$ be a flow on $X$, whose fixed point set $S_v$ has finitely
many connected components.\\
\indent In this paper, we always assume that $\mathcal{A}$ is an
isolated invariant set of $v$ containing $S_v$ and has finitely
many connected components. Denote its connected components
(including
$S_v$) by $A_i, i=1,2,\cdots,n$.\\
\indent \indent Let $A_i\in\mathcal{A}$ be a component with an
isolating neighborhood $U_i$, then for any given closed
neighborhood $V_i$ of $A_i$ in $U_i$, there is a flow neighborhood
(see \cite{FJ2} for definition) $A_i(r):=A_i(V_i,r)$. These are
compact
sets in $X$.\\
\indent Consider the flow $v$ restricted to the open space $X
\slash \cup^n_{i=1}A_i(r)$. Then all the trajectories can be
classified into three types:\\
\indent (1) trajectories with domain $(a, b),
-\infty<a<b<+\infty$;\\
\indent (2) trajectories with domain $(a, +\infty)$ or $(-\infty,
b), -\infty<a,b<+\infty$;\\
\indent (3) $(-\infty, +\infty)$.\\
\noindent We denote the sets of three type trajectories by
$\G^{\mathcal{A}}_1(r),\G^{\mathcal{A}}_2(r)$ and
$\G^{\mathcal{A}}_3(r)$, respectively. If $A_i$ is a point $p_i$
in $X$ for $i=1,\cdots,n$, then we denote $\G^{\mathcal{A}}_\g(r)$
by $\G_\g(r),\g=1,2,3$.

If $A\subset X$ is an open set such that $[\a]|_A=0$, then there
is a continuous function $\b$ such that $\a=\d\b$. Define
$I_\a(x,y):=\b(y)-\b(x)$, for $x,y\in A$. It is obvious that
$I_\a(x,y)$ depends only on $[\a]$.\\
\ \\
{\em {\bf Definition 2.1}$\;\;$ The flow $v$ defined on a compact
metric space $(X,d)$ is said to be a generalized $\a-$flow with
respect to an isolated invariant set
$\mathcal{A}=\{A_n,\cdots,A_1\}$,
 if there exist a continuous 1-cocycle $\a$, a small $r>0$ and a $T_0>0$ such that
for some $\r>0$ and $0\le\la<1$, the following conditions are satisfied:\\
\indent 1) $[\a]|_{A_i}=0,\;\;\max_{(x,y)\in A_i(r)\times
A_i(r)}|I_\a(x,y)|\le
\la\r,$ for $1\le i\le n $.\\
\indent 2) for any trajectory $\g\in \G^{\mathcal{A}}_1(r)$,
$$
\int_\g\a\ge\r
$$
\indent 3) if $\g(T_0)$ denotes any sub-trajectory of $\g\in
\G^{\mathcal{A}}_2(r)\cup \G^{\mathcal{A}}_3(r)$ with time
interval $T_0$, then
$$
\int_{\g(T_0)}\a\ge \r
$$
 \\
{\bf Remark 2.1}\quad In fact, we can take the constant $\r$ in
condition (3) to be different from the one in conditions (1) and
(2). This broader assumption will not change any proofs and
conclusions in \cite{FJ2}. Since $\a$ is a continuous cocycle, the
three conditions in Definition 2.1 also hold for $r'$ sufficiently
close to $r$. \\
 \\
\noindent {\bf Definition 2.2} $\;\;$ If the set $\mathcal{A}$ in
the definition of a generalized $\a$-flow $v$ consists only of
points, then $v$ is called an $\a$-flow. If $v$ is an $\a$-flow or
generalized $\a$-flow, we simply call it a flow carrying a cocycle
$\a$.
}  \\
 \\
\indent The definition of a continuous 1-cocycle in bounded
Alexander-Spanier cohomology theory and its integration along a
curve or a chain was already defined in \cite{Ch,FJ2} .\\
\indent Some interesting and important (generalized) $\a$-flows are given below, more examples can be found in \cite{FJ2}.\\
 \\
\noindent {\bf Example 2.1}$\;\;$ In classical Novikov theory, the
closed Morse 1-forms or the Bott-type closed 1-form were
considered frequently. Now in the present example, we consider the
flow $V$ generated by a closed 1-form $\o$ without any
non-degeneracy condition, even without a Bott-type condition. This
case will also introduce the Novikov type inequalities to be presented in the last section.\\
\indent Take the set $\mathcal{A}$ to be the zero locus of $\o$
and to have connected components $\{A_n,\cdots,A_1\}$.\\
\indent Firstly we show that any $A\in\mathcal{A}$ is an isolated
invariant set. Choose an open covering $\mathcal{U}=\{U_\a\}$ of
$A$ such that $\o=df_\a$ on $U_\a$. Let $U\subset\bigcup U_\a$ be
a tubular neighborhood of $A$ and $\pi:\;U\rightarrow A$ be the
projection. We want to show $I(U)=A$. Let $p\in I(U)$, then for
any $t_0, t_1, t_0<t_1$,
\begin{align*}
&\int_{t_0}^{t_1}|\o(p\cdot t)|^2 dt=\int_{[p\cdot t_0,p\cdot
t_1]}\o=\int_{[p\cdot t_0, \pi(p\cdot t_0)]}\o\\
&+\int_{[\pi(p\cdot t_0),\pi(p\cdot t_1)]}\o+\int_{[\pi(p\cdot
t_1), p\cdot t_1]}\o\\
&\le C
\end{align*}
We use the fact that the integration of $\o$ is only dependent on
the relative homotopy class of $[p\cdot t_0, p\cdot t_1]$.\\
\indent The above inequality implies that there exist sequences
$t_{0k}\rightarrow -\infty$ and $t_{1k}\rightarrow +\infty$ such
that $|\o(p\cdot t_{0k})|\rightarrow 0$ and $|\o(p\cdot
t_{1k})|\rightarrow 0$ as $k\rightarrow \infty$. Now the above
formula forces $|\o(p\cdot t)|=0,\forall t\in{\Bbb R}$. This
proves that $p\in A$.\\
\indent From the above conclusion, we can take a closed
neighborhood $V_i$ of $A_i$ contained in an isolated neighborhood
$ U_i$ s.t. $d(\pat V_i, \pat U_i)>r_0,
\min_{x\in\overline{U_i-V_i}}|\o(x)|\ge\ape_0$. Now take the flow
neighborhood $A_i(r):=A_i(B_{2r}(A_i),r)$ for small $r$. Since any
$\g\in\G^{\mathcal{A}}_1(r)$ has to go across $U_i-V_i$ for some
$i$, we have
\begin{align*}
\int_\g\o&=\int^b_a|\o(\g(t))|^2\;dt\\
&\ge \min_{1\le i\le n}(\min_{x\in\overline{U_i-V_i}}|\o(x)|\cdot
d(\pat V_i, \pat U_i))\ge r_0\ape_0
\end{align*}
Now if we take $r$ sufficiently small, then (1) and (2) holds.\\
\indent To prove (3) in the definition of flows carrying a
cocycle., we fix $r$ and then choose $T_0$ sufficiently large such
that $\forall \g\in\G^{\mathcal{A}}_2(r)\cup
\G^{\mathcal{A}}_3(r)$,
\begin{align*}
\int_{\g(T_0)}\o=\int |\o(\g(t))|^2\;dt\ge
\min_{x\in\overline{X-B_r(A_i)}}|\o(x)|^2 T_0.
\end{align*}
\indent As a direct corollary, the flow generated by a closed
1-form with Bott-type
non-degenerate zero sets is a generalized $\o$-flow w.r.t. the zero locus of $\o$.\\
 \\
{\bf Example 2.2}$\;\;$ A flow $v$ on a manifold $M$ is called a Morse-Smale flow (after \cite{Sm}) if it satisfies:\\
\indent (1) The chain recurrent set of $v$ consists of a finite
number of hyperbolic closed orbits and hyperbolic fixed points
.\\
\indent
(2) The unstable manifold of any closed orbit or fixed point has transversal intersection with the stable manifold of any closed orbit or fixed point.\\
\indent
Smale \cite{Sm} proved that such flows have ``global'' gradient-like structures and have a Morse decomposition which induces the Morse inequalities.\\
\indent However, in some cases, although the nonwandering set of
the flow contains only the hyperbolic periodic orbits and the
hyperbolic fixed points, the flow is not a Morse-Smale flow
because of the existence of ``cycles'' which consist of some
orbits ``connecting'' different invariant sets and form a closed
curve.   We can give a
 definition of such flows when restricted to the category of flows carrying a cocycle.\\
 \\
{\em {\bf $\a$-Morse-Smale flow}$\;\;$ Let $v$ be a generalized
$\a$-flow with respect to
 an isolated invariant set $\mathcal{A}=\{A_n,\cdots,A_1\}$. If $\mathcal{A}$ contains only the hyperbolic orbits
 or hyperbolic fixed points, then $v$ is called an $\a$-Morse-Smale flow.}\\
  \\
\indent Example 3.3.9 in \cite{FJ2} provides concrete $\a$-Morse-Smale flows.\\
 \\
\noindent {\bf Example 2.3 (Flows carrying a cohomology
class)}$\;\;$ We consider an extreme case namely that
$\mathcal{A}$ is empty in the definition of an $\a$-flow $v$. In
this case, the set $\G_1(r)\cup\G_2(r)=\emptyset$ and the only
condition that makes $v$ an $\a$-flow is that there exist
constants $\r>0$ and $T_0>0$ such that for any trajectory
$\g(T_0)$ with time interval $T_0$, we have
\begin{align*}
\int_{\g(T_0)}\a\ge \r
\end{align*}
\indent Now the $\a$-flow $v$ becomes a so called ``flow carrying
a cohomology class'' as introduced by R.C.Churchill \cite{Ch}. The
reason that the flow is called ``carrying a cohomology class'' is
that the above condition is independent of the choice of the
representative in the cohomology class $[\a]$. In fact, if
$\a_1\in[\a]$ is another cocycle, then there exists a coboundary
$\d\b$, such that
\begin{align*}
\a_1-\a\simeq \d\b
\end{align*}
and so
\begin{align*}
\int_{\g(kT_0)}\a_1&=\int_{\g(kT_0)}\a+\b(e(\g(kT_0)))-\b(s(\g(kT_0)))\\
&\ge k\r-2M_\b
\end{align*}
where $s(\g)$ and $e(\g)$ are the start point and the end point of
the trajectory $\g$, and $M_\b$ is the bound of $\b$. Hence if we
choose $k>[\frac{2M_\b+\r}{\r}]+1$, then we have
\begin{align*}
\int_{\g(kT_0)}\a_1\ge \r.
\end{align*}
\indent
The existence of a flow carrying a cohomology class in a manifold will induce the vanishing of the Novikov numbers.
This result will be given in the last section.\\
 \\
{\bf  3. $\a$-flows and gradient-like flows}\\
 \\
After defining the $\a$-flows, it is natural to find the relation
between $\a$-flows and the well-known gradient-like flows. It
seems that an $\a$-flow for $\a$ being a nontrivial cocycle should
be a nongradient-like flow. However the following simple example
shows that this is not true.\\
 \\
{\bf Example 3.1}\quad Let $S^1$ be the unit circle with the
standard metric. Let $\phi_1$ and $\phi_2$ be two strictly
monotonically increasing functions in the interval $[0,\pi]$
satisfying the following conditions
\begin{align*}
&\phi_1(0)=0,\phi'_1(0)=0;\;\phi_1(\pi)=3\pi,\phi'_1(\pi)=0\\
&\phi_2(0)=0,\phi'_2(0)=0;\;\phi_2(\pi)=\pi,\phi'_2(\pi)=0.
\end{align*}
Define a smooth circle-valued function
$$
f(x)=\left\{\begin{array}{ll}
e^{i\phi_1(x)}&0\le x\le \pi\\
e^{i\phi_2(2\pi-x)}& \pi\le x\le 2\pi.\end{array}\right.
$$
Consequently, $f(x)$ induces a cocycle $\alpha
(x,y):=J(f)(x,y):S^1\times S^1\to {\Bbb R}$ such that if $y$ lies
in a small neighborhood of $x$, then
$$
\int_{[x,y]}\alpha=\text{arg}f(y)-\text{arg}f(x).
$$
Here the map $J$ is defined in proposition 3.1.1 of \cite{FJ2}.
This cocycle is a nontrivial cocycle. Since if we let $\g$ be the
oriented curve starting from the point $\ta=0$ and then going
around the circle in the anticlockwise direction, then for any
integer $l$, we have
$$
\int_{l\g}\a=2\pi l.
$$
\indent Now we consider the gradient-like flow $v$ on $S^1$ that
has two fixed points at $\ta=0,\pi$, and flows from the point
$\ta=0$ to the point $\ta=\pi$. Then it is easy to check that $v$
carries the cocycle $\a$ with parameter $\la=0,\r=1$ and
$\G_2(r)\cup \G_3(r)=\emptyset$ in the definition 2.1.

Thus, this example demonstrates that a gradient-like flow can
carry a nontrivial cocycle.
 \\
{\bf Remark 3.1}\quad It is easy to see that for any $x\in S^1$,
and $t>0$, $\int_{[x,x\cdot t]} \alpha>0$. Hence the pull-back
form $f^*(d\ta)$ is a Lyapunov 1-form for the pair $(-v,
\ta=0,\pi)$, which is defined in \cite{Fa3, Fa4}. Though the
smooth closed 1-form $f^*(d\ta)$
is a nontrivial smooth closed 1-form, $-v$ is a gradient-like flow. \\

We will give a series of theorems and examples to formulate the
relation between $\a$-flows and gradient-like flows.
The following theorem is already given in \cite{FJ2}.\\
 \\
{\em {\bf Theorem 3.1}$\;\;$ Let $v$ be an $\a$-flow on the
compact metric space $(X,d)$. If $\a$ is a trivial cocycle, then
$v$ is a gradient-like flow. Conversely, if $v$ is a gradient-like
flow, then $v$ is a $\d g$-flow for some continuous function
$g$.\\}
 \\
We will give necessary and sufficient conditions to distinguish
$\a$-flows and gradient-like flows. Firstly we study the local
structure of a gradient-like flow. \\
 \\
{\bf Local structure of gradient-like flows}\quad Let $p$ be a
fixed point of a gradient-like flow $v$. Then there
 exists an associated Lyapunov function $g(x)$ on a small neighborhood $\bar{B}_r(p)$.\\
\indent Let $0<s<\frac{r}{2}$ and let $B_{\frac{r}{2}}(p)$ be a
closed ball centered at $p$ with radius $\frac{r}{2}$. Define two
sets on the sphere $\pat B_r(p)$,
$$
\aligned
&B^+_{r,s}(p)=\{x\in\partial B_r(p)|[x,x\cdot t]\cap \pat B_s(p)\neq\emptyset\;\hbox{for some}\; t>0\}\\
&B^-_{r,s}(p)=\{x\in\partial B_r(p)|[x,x\cdot t]\cap \pat
B_s(p)\neq\emptyset\;\hbox{for some}\; t<0\}
\endaligned
$$
\indent For any point $x\in B^+_{r,s}(p)$, let $t_x$ denotes the
reach time of the trajectory $[x, x\cdot t](t>0)$ to the sphere
$\partial B_{\frac{r}{2}}(p)$. We have the following lemma.\\
 \\
{\em {\bf Lemma 3.1}\quad (1) $B^\pm_{r,s}(p)$ are closed sets on
$\partial B_r(p)$.\\
(2) $t_x$ is a lower semicontinuous function on $B^+_{r,s}(p)$.}\\
 \\
 Proof. $B^\pm_{r,s}(p)$ are actually the boundaries of the flow neighborhoods
  $A^\pm (B_r, s)$ defined in \cite{FJ2}. Hence they are closed.
  For convenience, we give a proof here.
  To prove (1), we need only to prove that the set $\partial
 B_r(p)-B^+_{r,s}(p)$ is open in $\partial B_r(p)$ and will
 finally drop into another fixed point $q$. In terms of the choice
 of $g(x)$, we can choose $a<b$ satisfying
 \begin{equation*}\label{$3.1$}
 B_r(p)\subset g^{-1}([b,+\infty]);\;B_r(q)\subset
 g^{-1}((-\infty,a])\tag{$3.1$}
\end{equation*}
Since $x\cdot t(t>0)$ drops into $q$, there exists a $T_x>0$ such
that $x\cdot T_x$ is in the interior of $B_r(q)$. Therefore there
is a closed ball $D_\delta(x)$ such that for any $y\in
D_\d(x)\cap\partial B_r(p), y\cdot T_x$ is in the interior of
$B_r(q)$. We can choose $\d$ small enough such that the set
$D_\d(x)\cdot [0,T_x]\cap \partial B_s(p)=\emptyset$. Since
$y\cdot T_x\in g^{-1}((-\infty,a]),y\cdot [T_x,+\infty)\in
g^{-1}((-\infty,a])$. Therefore $D_\d(x)\cdot [0,+\infty]\cap
\partial B_s(p)=\emptyset$ in view of (\ref{$3.1$}). This shows that
$\partial B_r(p)-B^+_{r,s}(p)$ is open in $\partial B_r(p)$.

To prove (2), let $x\in B^+_{r,s}(p)$ and let $t_x<+\infty$ be the
arrival time. We need only to prove that for any small
$\epsilon>0$, there exists a neighborhood $D_\delta(x)$ of $x$
such that for any $y\in D_\delta(x)\cap B^+_{r,s}(p)$,
\begin{equation}\label{cont1}
t_y\ge t_x-\epsilon\tag{3.2}
\end{equation}
Since $t_x$ is the arrival time of the trajectory $x\cdot t(t>0)$
to $\partial B_{\frac{r}{2}}(p)$, $x\cdot [0,t_x-\epsilon]$ has a
positive distance from $\partial B_{\frac{r}{2}}(p)$. By the
continuity of the flow $v$, there is a neighborhood $D_x$ such
that $\bar{D}_x\cdot [0,t_x-\epsilon]$ has a positive distance
from $\partial B_{\frac{r}{2}}(p)$ as well. Hence for $y\in
D_x\cap
B_{r,s}^+(p)$, (\ref{cont1}) is true. \qed\\

 Since the Lyapunov function $g(x)$ is strictly decreased along any
 nonconstant trajectory, there exists a $\d_x>0$ such that
$g(x)-g(x\cdot t_x)=2\d_x$. By lemma 3.1, (2) and the property of
$g(x)$, the function $g(y)-g(y\cdot t_y)$ is lower semicontinuous
with respect to $y$ on $B^+_{r,s}(p)$. Hence for each $x\in
B^+_{r,s}(p)$, there exists a neighborhood $U_x$ in $
B^+_{r,s}(p)$ such that $\forall y\in U_x,\,g(y)-g(y\cdot
t_y)>\d_x$. This implies that $\d^+_{r,s}(p):=\min_{x\in
B^+_{r,s}(p)}{\d_x}$ is positive and $\forall x\in
B^+_{r,s}(p),\,g(x)-g(x\cdot t_x)>\d^+_{r,s}(p)$. In the same way,
we can obtain a $\d^-_{r,s}(p)>0$ such that
$\forall x\in B^-_{r,s}(p),$ there is $g(x\cdot(-t_x))-g(x)>\d^-_{r,s}(p)$.\\

Let $0<s_1<s_2<\frac{r}{2}$, then we have the following
conclusions:
$$
\aligned
&(1)\;\;B^+_{r,s_1}(p)\subset B^+_{r,s_2}(p);B^-_{r,s_1}(p)\subset B^-_{r,s_2}(p)\\
&(2)\;\;\d^\pm_{r,s_2}(p)\le \d^\pm_{r,s_1}(p)\\
&(3)\;\;\mbox{There exists $s_0>0$, such that for any $0<s<s_0$, }
B^+_{r,s}(p)\cap B^-_{r,s}(p)=\emptyset
\endaligned
$$
The first two conclusions are obvious. For (3), we can choose
$s_0>0$ satisfying $\osc_{x\in
B_{s_0}(p)}g(x)<\min{\d^\pm_{r,s_0}(p)}$. If there is a $0<s<s_0$
such that $ B^+_{r,s}(p)\cap B^-_{r,s}(p)\neq\emptyset$, then this
means that there is a trajectory in $\bar{B}_r(p)$ which has non
empty intersection with $\partial B_r(p)$ and $\partial B_s(p)$.
But this is absurd, because $g(x)$ is strictly decreacing along
any trajectory. After travelling from
 $\partial B_s(p)$ to $\partial B_r(p)$, then back to $\partial B_s(p)$, the value
 of $g(x)$ will decrease at least $\d^+_{r,s}(p)+\d^-_{r,s}(p)$ which contradicts
  the fact that $\osc_{x\in B_s(p)}g(x)\le\osc_{x\in B_{s_0}(p)}g(x)<\min{\d^\pm_{r,s}(p)}$.
   Hence for $0<s<s_0$, (3) holds.\\
 \\
{\bf Integration of cocycles in a gradient-like flow}\quad In this
part we assume that $v$ is a gradient-like flow with finitely many
fixed points $\{p_j;j=1,2,\cdots,n\}$. Let $\a$ be a cocycle on
$X$. Then there is an associated $r$-covering
$\mathcal{U}(r)$(covering consisting of radius $r$-balls) such
that when restricted to the closure of each open ball
 in $\mathcal{U}(r)$, $\a$ is a coboundary $\d\b$, where $\b$ is a bounded function.
 Let $r$ be small enough such that the Lyapunov function $g(x)$ corresponding to the flow $v$
has mutually non-intersecting images $g(B_r(p_j))$ for $j=1,2,\cdots,n$.\\
\indent Define $\d_0(r)=\min_{1\le j\le
n}\d^\pm_{r,\frac{r}{4}}(p_j)$. Since $g(x)$ is uniformly
continuous, there exists $s_0\le \frac{r}{4}$ satisfying
$\osc_{B_{s_0}}{g(x)}<\frac{1}{3}\d_0(r)$. By our choice and in
view of the analysis of the local structure of gradient-like
flows, it is easy to see that the following inequality holds:
\begin{equation}\label{3.3}
3\max_x\osc_{y\in
B_{s_0}(x)}g(y)<\d_0(r)<\d^\pm_{r,\frac{r}{4}}(p_j)\le
\d^\pm_{r,s_0}(p_j)\tag{3.3}
\end{equation}
for $ j=1,2,\cdots,n$.\\
 \\
{\em {\bf Theorem 3.2}$\;\;$ Let $v$ be a gradient-like flow with
finitely many fixed points. $\a$ is a cocycle with an associated
$r$-covering. Let $\mathcal{U}(s)$ be an $s$-covering with
$s<s_0<\frac{r}{4}$ for some $s_0$ satisfying (\ref{3.3}). Then
there exist $M>0$ and $T(s_0)>0$ depending on $r,n,\a$ but not on
$s$ such that for any $(\mathcal{U}(s), T(s_0))$-chain
$\tilde{\g}$, we have
\begin{equation}\label{3.4}
\Big|\int_{\tilde{\g}}\a\Big|\le M \tag{3.4}
\end{equation}
}\\

Recall that a $(\mathcal{U}(s),T)$-chain from $x$ to $x'$ is a
sequence $\{x=x_1,\cdots,x_{n+1}=x'|t_1,\cdots,t_n\}$ such that
$t_i\ge T$ and each pair $(x_i\cdot t_i, x_{i+1})(i=1,\cdots,n)$
belongs to a ball with radius $s$ in the covering $\mathcal{U}$.\\
 \\
Proof$\;\;$ Firstly we will prove that given any $ s\le r$, for any trajectory $\g\in\G_1(s)$ with domain $[a,b]$, there exists a constant $T(s)$ satisfying $b-a\le T(s).$\\
\indent
Let $x\in\partial B_s(p_j)$, then the trajectory $x\cdot t$ or $x\cdot(-t)$ for $t\ge 0$ will flow into some different fixed point of $v$. (There does not exist any trajectory joining one point to itself.) Without loss of generality, we assume that $(x\cdot t)(t\ge 0)$ flows into the point $p_k$ for $k\neq j$. (It is possible that the trajectory may pass through some ball $B_s(p_l)$, but this does not influence the result below.) Hence there exists a $T_x>0$ such that $x\cdot T_x\cap B_{\frac{s}{2}}(p_k)\neq\emptyset$. By the continuity of the flow, there exists a small closed ball $D_{s_x}(x)$ which satisfies that for any point $y\in D_{s_x}(x)\cap\pat B_s(p_j)$, the trajectory $[y, y\cdot T_x]\cap\pat B_s(p_k)\neq\emptyset$. This shows that for any trajectories starting from $D_{s_x}(x)\cap \pat B_s(p_j)$ their time intervals are not greater than $T_x$, and covering the sphere $\pat B_s(p_j)$ by finitely such small closed balls, then it is easy to  see that any trajectory in $\G_1(s)$ starting from $\pat B_s(p_j)$ has time interval not greater than a constant $ T_j(s)$. Let $T(s)=\max_{1\le j\le n}{T_j(s)}$, then all the trajectories in $\G_1(s)$ have time interval not greater than $T(s)$.\\
\begin{figure}
\begin{center}
\includegraphics{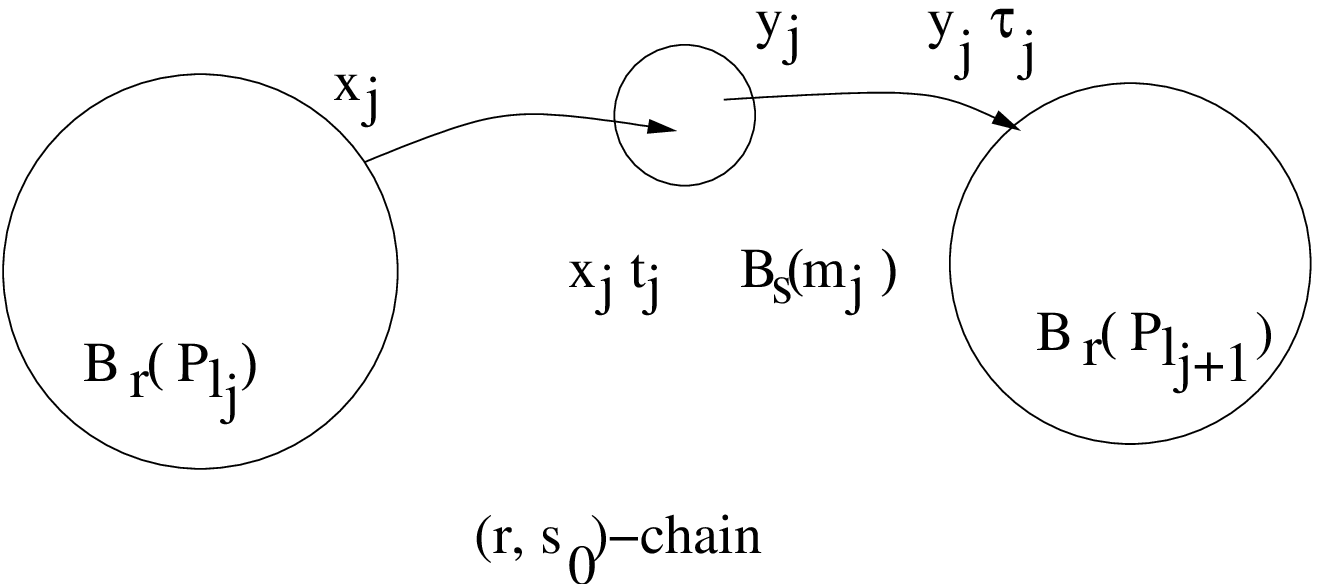}
\end{center}
\end{figure}
\indent
Let $\g=\{x_0,y_0,x_1,y_1,\cdots,x_{k-1},y_{k-1},x_k|t_0,\tau_0,\cdots,t_{k-1},\tau_{k-1}\}$ be a chain. It
is called an $(r,s_0)$-chain if it satisfies the following conditions:\\
\begin{description}
\item[(1)] $x_j\in\pat B_r(p_{l_j})$ for $j=0,1,\cdots,k-1$ and
$x_k\in B_r(p_{l_k})$;

\item[(2)] $y_j\cdot\tau_j\in\pat B_r(p_{l_{j+1}})$ for
$j=0,1,\cdots,k-1$;

\item[(3)] the point pair $(x_j\cdot t_j,y_j)\in B_s(m_j)$, for
$j=0,1,\cdots,k-1$, where those $B_s(m_j)'s$ are elements in the
$s$-covering $\mathcal{U}(s)$ of $X$.

\item[(4)] $B_s(p_{l_j}),j=0,1,\cdots,k$ are $k$ different balls,
hence $k\le n-1$.

\item[(5)] $t_j,\tau_j\le T(s_0)$, where $T(s_0)$ is an upper
bound for the time interval of all the trajectories in
$\G_1(s_0)$.

\end{description}
Now compute the integral of $\a$ along the $(r,s_0)$-chain $\g$.\\
\begin{equation}\label{3.5}
\aligned
\Big|\int_{\bar{\g}}\a\Big|&\le \sum^{k-1}_{j=0}\Big|\int_{[x_j,x_j\cdot t_j]}\a\Big|+\sum^{k-1}_{j=0}\Big|\int_{[y_j,y_j\cdot\tau_j]}\a\Big|\notag\\
&+\sum^{k-1}_{j=0}\Big|\int_{[x_j\cdot t_j,y_j]}\a\Big|+\sum^{k-1}_{j=0}\osc_{x\in B_r(p_{l_{j+1}})}{\b_{p_{l_{j+1}}}}\notag\\
&=I+I\!\!I+I\!\!I\!\!I+I\!\!V
\endaligned\tag{3.5}
\end{equation}
Since $0<t_j,\tau_j\le T(s_0)$, applying corollary 3.2.3 in
\cite{FJ2},
$$
I+I\!\!I\le 2kC_1
$$
where $C_1$ is a constant depending only on $\max{|\dot{v}(t)|},\a$ and $r$.\\
\indent Connect the point pair $(x_j\cdot t_j,y_j)$ with a line
segment, then its length is at most $2s$. Using proposition 3.2.2
in \cite{FJ2}, we have
$$
I\!\!I\!\!I\le kC_2
$$
Here $C_2$ depends only on $s_0,\a$ and $r$.\\
\indent For $I\!\!V$,
$$
I\!\!V\le 2kM_\a(r)
$$
Combining the above estimates, we have
\begin{equation}\label{3.6}
|\int_{\bar{\g}}\a|\le kC\tag{3.6}
\end{equation}
where $C$ depends on $r,s_0,\a$ and $\max{|\dot{v}(t)|}$.\\
\indent
Now the proof of the proposition is changed to the problem to reduce each
$(\mathcal{U}(s),T(s_0))$-chain to an $(r,s_0)$-chain, while keeping the integral of $\a$ invariant.\\
\indent
Let $\tilde{\g}=\{x_0,x_1,\cdots,x_k|t_0,\cdots,t_{k-1}\}$ be a $(\mathcal{U}(s),T(s_0))$-chain.
Consider the following cases.\\
\begin{figure}
\begin{center}
\includegraphics{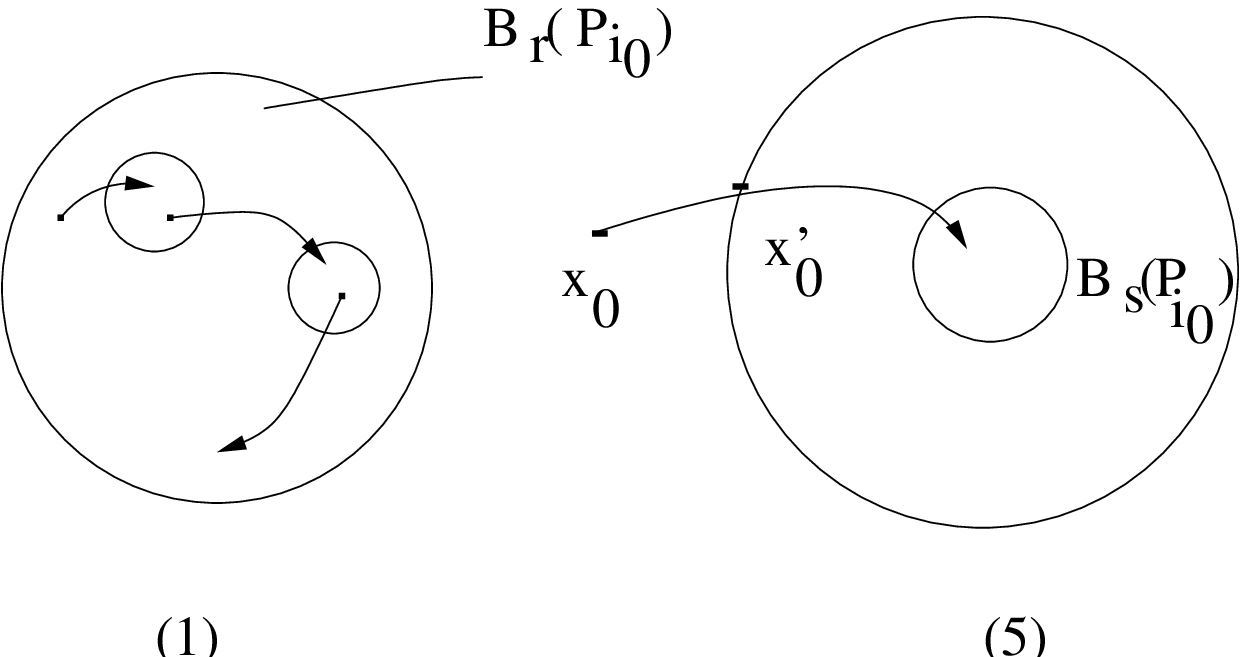}
\end{center}
\end{figure}
\indent (1) If the chain $\tilde{\g}$ is contained in a ball
$\bar{B}_r(p_{i_0})$ for $0\le i_0\le n$, then
\begin{align*}
&|\int_{\tilde{\g}}\a|=|\b_{p_{i_0}}(x_k)-\b_{p_{i_0}}(x_0)|\\
& \le\osc_{x\in B_r(p_{i_0})}{\b_{p_{i_o}}(x)}\le 2M_\a(r)
\end{align*}
\begin{figure}
\begin{center}
\includegraphics{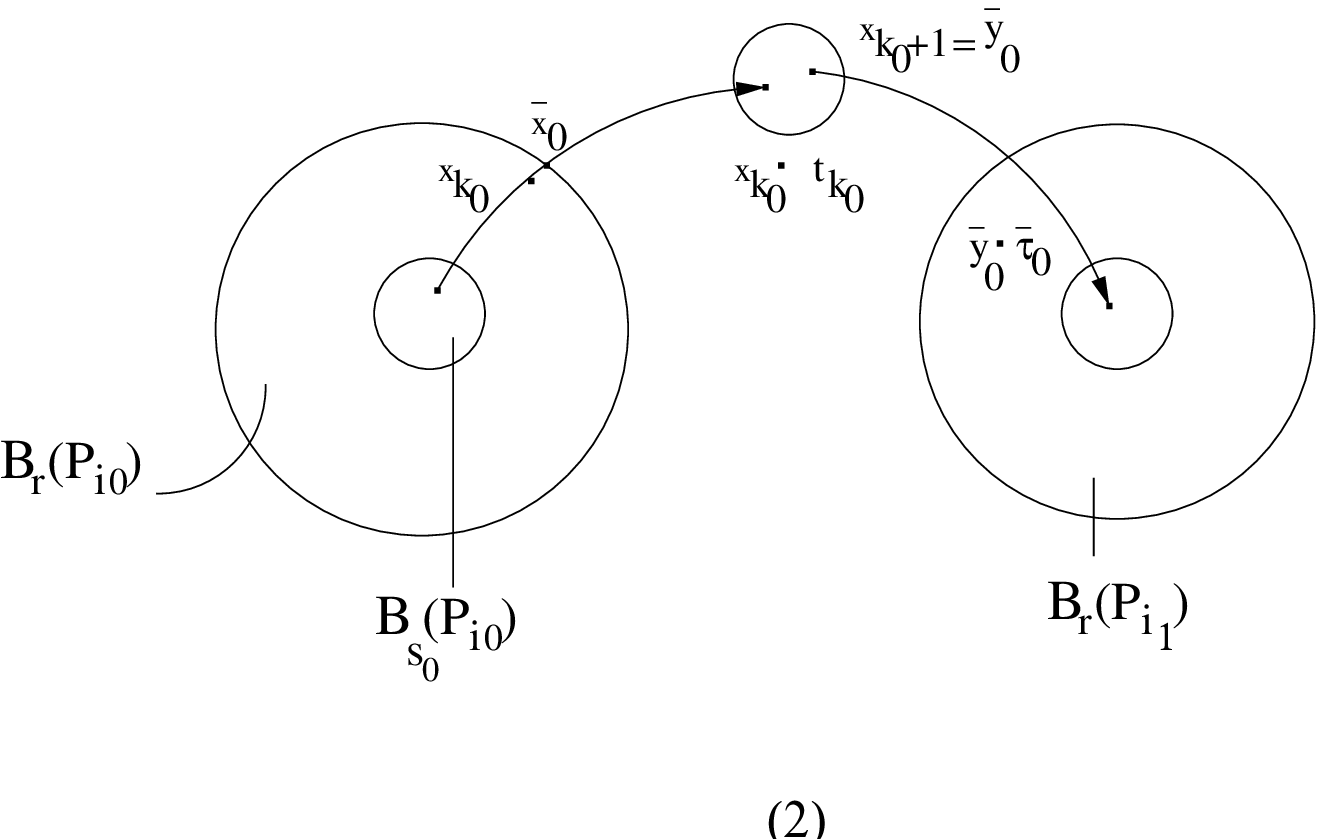}
\end{center}
\end{figure}
\indent (2) There is $1\le k_0\le k$ such that the subchain
$\tilde{\g}_{0k_0}=
\{x_0,x_1,\cdots,x_{k_0}|t_0,\cdots,\\t_{k_0-1}\}$ is contained in
$B_r(p_{i_0})$
 but the point $x_{k_0}\cdot t_{k_0}$ is not in $\cup^n_{j=0}\bar{B}_r(p_j)$. Hence
  there is a first intersection point $\bar{x}_0$ of $[x_{k_0},x_{k_0}\cdot t_{k_0}]$
   and $\pat B_r(p_{i_0})$. Let $\bar{t}_0$ satisfy $\bar{x}_0\cdot\bar{t}_0=x_{k_0}
   \cdot t_{k_0}$ and let $\bar{y}_0=x_{k_0+1}$. Denote the ball containing the point
   pair $(\bar{x}_0\cdot\bar{t}_0,\bar{y}_0)$ by $B_s(m_0)$. Since $t_{k_0}\ge T(s_0)$
    and $\tilde{\g}_{0k_0}\subset B_r(p_{i_0})$, this implies that $x_{k_0}\in B_{s_0}(p_{i_0})$,
     hence the trajectory $[\bar{x}_0,\bar{x}_0\cdot\bar{t}_0]$ does not intersect with
     $B_{s_0}(p_{i_0})$. Otherwise,$[x_{k_0},x_{k_0}\cdot t_{k_0}]$ starts from $B_{s_0}
     (p_{i_0})$, intersects $\pat B_r(p_{i_0})$, then goes back to  $B_{s_0}(p_{i_0})$.
     This conclusion contradicts (\ref{3.3}) when we check the change of the Lyapounov
     function $g(x)$ along $[x_{k_0},x_{k_0}\cdot t_{k_0}]$. Therefore $[\bar{x}_0,
     \bar{x}_0\cdot\bar{t}_0]$ is part of a trajectory in $\G_1(s_0)$, and we have
     $\bar{t}_0\le T(s_0)$. Since $\bar{x}_0\cdot\bar{t_0}=x_{k_0}\cdot t_{k_0}$ is not
      in $\bar{B}_r(p_{i_0})$, by the choice of $s$, we know that $\bar{y}_0:=x_{k_0+1}
      \notin B_{s_0}(p_{i_0})$. Hence $\bar{y}_0\cdot t_{k_0+1}$ meets firstly $B_{s_0}
      (p_{l_1})$ for some fixed point $p_{l_1}$. Let $\bar{y}_0\cdot \bar{\tau}_0$ be the
       last intersection point of $[\bar{y}_0,\bar{y}_0\cdot t_{k_0+1}]$ with $\pat B_r(p_{i_1})$.
        In the same way, we can prove that $\bar{\tau}_0\le T(s_0)$. Joining
        $x_0$ and $\bar{x}_0$ by a curve $l_{x_0\bar{x}_0}$, then the subchain $\tilde{\g}_{0k_0}$
        is reduced to a curve $l_{x_0\bar{x}_0}$ combined with a chain
        $\{\bar{x}_0,\bar{y}_0|\bar{t}_0,\bar{\tau}_0\}$ and it is clear that
        the reduction keeps the integral of $\a$ invariant.\\
\begin{figure}
\begin{center}
\includegraphics{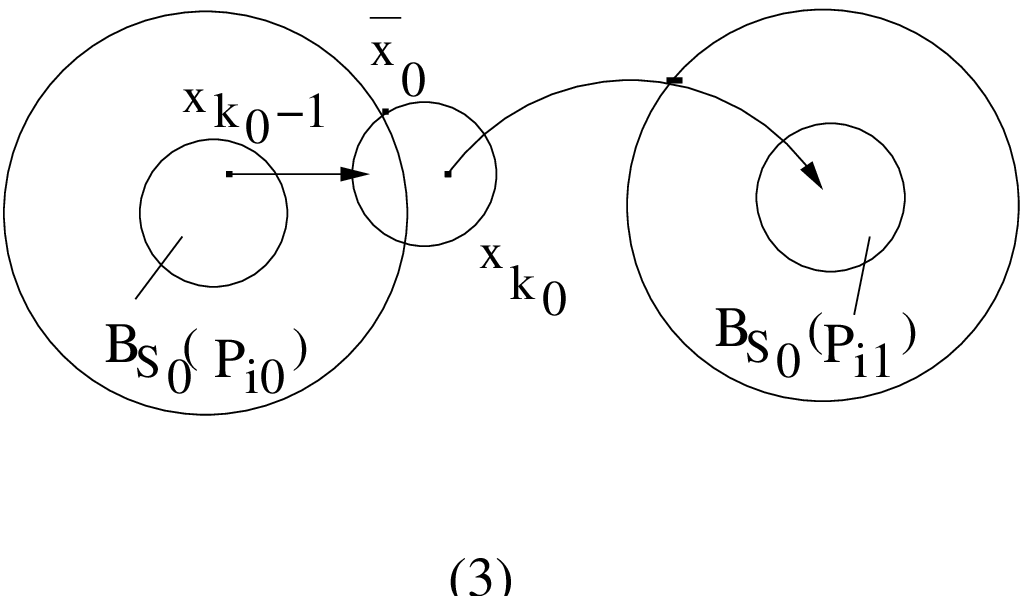}
\end{center}
\end{figure}
\indent (3) There is a $2\le k_0\le k$ such that the subchain
$\tilde{\g}_{0,k_0-1}
=\{x_0,x_1,\cdots,x_{k_0-1}|t_0,\\t_1,\cdots,t_{k_0-2}\}$ and the
trajectory
 $[x_{k_0-1},x_{k_0-1}\cdot t_{k_0-1}]$ is contained in $B_r(p_{i_0})$ but
 $x_{k_0}\notin B_r(p_{i_0})$. Then we take a point $\bar{x}_0\in\pat B_r(p_{i_0})
 \cap B_s(m_0)$, where $B_s(m_0)$ is the ball containing $(x_{k_0-1}\cdot t_{k_0-1}
 ,x_{k_0})$. Since $t_{k_0}\ge T(s_0)$, the trajectory $[x_{k_0},x_{k_0}\cdot t_{k_0}]$
  will meet $B_{s_0}(p_{i_1})$. We denote $x_{k_0}$ by $\bar{y}_0$
  and represent the last intersection point of $[x_{k_0},x_{k_0}\cdot t_{k_0}]$
  with $\pat B_r(p_{i_1})$ by $\bar{y}_0\bar{\tau}_0$ for some $\bar{\tau}_0>0$.
  Due to the same argument as in (2), we know that $\bar{\tau}_0\le T(s_0)$.
  So $\tilde{\g}_{0k_0}$ is reduced to the combination of a curve $l_{x_0\bar{x}_0}$
  and a chain $\{\bar{x}_0,\bar{y}_0|0,\bar{\tau}_0\}$. It is clear that such
  inductions also keep the integral of $\a$ invariant. Note that $x_{k_0-1}$
  must be in $B_{s_0}(p_{i_0})$.\\
\begin{figure}
\begin{center}
\includegraphics{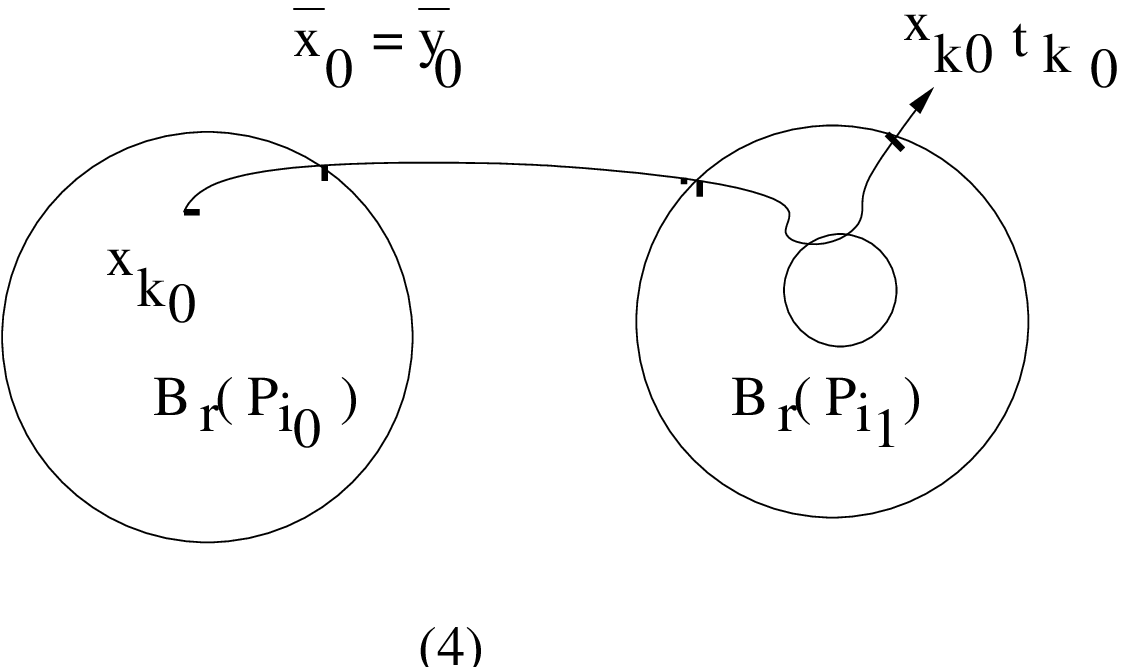}
\end{center}
\end{figure}
\indent (4) There is a subchain
$\tilde{\g}_{0k_0}=\{x_0,x_1,\cdots,x_{k_0}|t_0,\cdots,
t_{k_0-1}\}$ contained in $B_r(p_{i_0})$ and the trajectory $
[x_{k_0},x_{k_0}\cdot t_{k_0}]$ goes through $B_{s_0}(p_{i_1})$.
Let $\bar{x}_0=\bar{y}_0$ be the first intersection point of
$[x_{k_0},x_{k_0}\cdot t_{k_0}]$ with $\pat B_r(p_{i_1})$. Also we
can prove that $\bar{\tau}_0\le T(s_0)$ as in (2) and (3). Then
the chain $\tilde{\g}_{0k_0}$ followed by
$[x_{k_0},\bar{y}_0\cdot\bar{\tau}_0]$ is reduced to the
combination of $l_{x_0\bar{x_0}}$ with the
chain $\{\bar{x}_0,\bar{y}_0|0,\bar{\tau}_0\}$, while the integral of $\a$ is invariant under change.\\
\indent (5) If $x_0\notin\cup^n_{j=1}B_r(p_j)$, then $x_0\cdot
t_0$ will firstly meet some $B_s(p_{i_0})$. Let $x'_0=x_0\cdot
t'_0$ be the last intersection point
of $[x_0,x_0\cdot t_0]$ with $\pat B_r(p_{i_0})$. In the same way, we can prove $t'_0\le T(s_0)$.\\
\indent Now using the above steps (2)-(5) repeatedly, any
$(\mathcal{U}(s),T(s_0))$- chain
$\tilde{\g}=\{x_0,x_1,\cdots,x_k|t_0,\cdots,t_{k-1}\}$ can be
reduced to the combination of a curve $l_{x_0\bar{x}_0}$ ( or
$[x_0,x_0\cdot t_0']$ in case (5)) with an $(r,s_0)$-chain
$\bar{\g}=\{\bar{x}_0,\bar{y}_0,\cdots,\bar{x}_l|\bar{t}_0,
\bar{\tau}_0,\cdots,\bar{\tau}_{l-1}\}$ if we can prove that (4)
in the definition
of a $(r,s_0)$-chain holds.\\
\indent
If (4) is not true, then there exists an $(r,s_0)$-cycle $\bar{\g}_c=\{\bar{x}_{i_0},
\bar{y}_{i_0},\cdots,\bar{y}_{i_l},\bar{x}_{i_0}|\bar{t}_{i_0},\\ \bar{\tau}_{i_0},
\cdots,\bar{\tau}_{i_l}\}$ where the point pair $(\bar{y}_{i_l}\cdot\bar{\tau}_{i_l},
\bar{x}_{i_0})\in\pat B_r(p_{i_0})$. Now extend the trajectory $\bar{y}_{i_j}\cdot t(t\ge 0)$
 for $j=0,\cdots,l.$ Since $\bar{\tau}_c$ is obtained by reducing the $(\mathcal{U}(s),T(s_0))$
 -chain, hence $\bar{y}_{i_j}\cdot t$ will go through $B_{s_0}(p_{i_{j+1}})$.
 We let $\bar{y}_{i_j}\cdot\tau_{i_j}(\tau_{i_j}>\bar{\tau}_{i_j})$ be a point
 in $B_{s_0}(p_{i_{j+1}})$. Take the inverse process corresponding to (2)-(4),
 then there is a trajectory starting from some point $x_{i_j-1}\in B_{s_0}(p_{l_j-1})$
  to some point $x_{i_j-1}\cdot t_{i_j-1}$ in the $s$-ball containing
  $\bar{x}_{i_{j-1}}\cdot t_{i_{j-1}}$ and $\bar{y}_{i_j}$. Therefore from the
  cycle $\bar{\g}_c$ we get a $(\mathcal{U}(s),T(s_0))$-chain
  $\tilde{\g}_c=\{x_{i_0},\bar{y}_{i_0},\cdots,\bar{y}_{i_l},x_{i_0}|t_{i_0},
  \tau_{i_0},\cdots,t_{i_l},\tau_{i_l}\}$.\\
\indent Consider the change of the function $g(x)$ along
$\tilde{\g}_c$. On one hand, since the two ends of $\tilde{\g}_c$
are $x_{i_0}$ and $\bar{y}_{i_l}\cdot\tau_{i_l}$, by (\ref{3.3})
\begin{align*}
\Big|\int_{\tilde{\g}_c}\d g\Big|\le \osc_{x\in
B_{s_0}(p_{i_0})}g(x)<\frac{1}{3}\d_0(r)
\end{align*}
On the other hand, each $\bar{y}_{i_j}\cdot t
(0<t\le\tau_{i_j})(j=0,\cdots,l)$ goes across
$B_r(p_{i_j})\backslash B_{\frac{r}{4}}(p_{i_j})$, we have
\begin{align*}
\int_{\tilde{\g}_c}\d g\ge&
(l+1)\d^+_{r,\frac{r}{4}}-2(l+1)\max_{x}{\osc_{y\in B_s(x)}g(y)}
>&\d_0(r)-\frac{2}{3}\d_0(r)=\frac{1}{3}\d_0(r)
\end{align*}
This is absurd. The contradiction shows that the chain $\bar{\g}$
we get from a $(\mathcal{U}(s),T(s_0))$-chain $\tilde{\g}$ is
indeed an $(r,s_0)$-chain. Let
$\bar{\g}=\{\bar{x}_0,\bar{y}_0,\cdots,\bar{x}_l|\bar{t}_0,\bar{\tau}_0,\cdots,\bar{\tau}_{l-1}\}$.
Since $v$ is gradient-like, the index $l$ is less than $n$.
Therefore applying the estimate (\ref{3.6}) and the fact that the
reduction from $\tilde{\g}$ to $\bar{\g}$ keeps the integral of
$\a$, we have
$$
\Big|\int_{\tilde{\g}}\a\Big|=\Big|\int_{\bar{\g}}\a\Big|\le nC.
$$
where $C$ depends on $r, s_0, \a $ and $\max{|\dot{v}(t)|}$. Theorem 3.2 now is proved. \qed\\
 \\
{\em {\bf Theorem 3.3}$\;\;$ Let $v$ be an $\a$-flow on the
compact metric space $(X,d)$. Let $\{s_i\}\to 0$ as $i\to\infty$.
If $\a$ is a nontrivial cocycle and for any $M>0$ and $T>0$, there
is a $(\mathcal{U}(s_i),T)$-chain $\tilde{\g}_i$ for each $s_i$
such that
\begin{align}\label{3.7}
\left|\int_{\tilde{\g}_i}\a\right|\ge M\tag{3.7},
\end{align}
then $v$ is not a gradient-like flow. Furthermore, if it is known
that for any $M>0$, there is a trajectory $\tilde{\g}$ satisfying
(\ref{3.7}) or there is an oriented cycle $\tilde{\g}_0$
consisting of some orbits joining fixed points where the direction
is determined by the forward direction of the flow $v$ such that
\begin{align}\label{3.8}
\left|\int_{\tilde{\g}}\a\right|\ge 1\tag{3.8},
\end{align}
then $v$ is not a gradient-like flow.\\
}
 \\
Proof.\quad The first conclusion is a direct corollary of Theorem
3.2 and the third one is obvious. We only consider the second
case. It is easy to see that if the time intervals of the
trajectories have an upper bound, then for any cocycles on $X$ the
absolute value of the integral, $|\int_{\tilde{\g}}\a|$ has a
uniform bound. Therefore the trajectory $\tilde{\g}$ with respect
to the arbitrary large $M$ has arbitrary large time interval and
it is the chain needed for the hypothesis in the first
conclusion.\qed\\
 \\
{\em {\bf Theorem 3.4}$\;\;$ Let $v$ be an $\a$-flow with a
nontrivial cocycle $\a$ on a compact metric space $(X,d)$. If $v$
is not a gradient-like flow, then for any small sequence $\{s_i\}$
which tends to zero as $i\to\infty$, and for any $M>0, T>0$ there
is a $(\mathcal{U}(s_i), T)$-chain $\tilde{\g}_i$ for each $s_i$
such that
\begin{align}\label{3.9}
\left|\int_{\tilde{\g}}\a\right|\ge M\tag{3.9}
\end{align}
}\\
 \\
Proof\quad If $v$ is not a gradient-like flow, then there is a non
fixed point $x_0$ in the chain recurrent set of $v$. Therefore for
any $s$ and $T>0$ there is a $(\mathcal{U}(s),T)$-chain
$\tilde{\g}_i=\{x_0,\cdots, x_k=x_0|t_0,\cdots, t_{k-1}\}$.\\
\indent Now we consider two possibilities:\\
\indent (1). If there exists $x_i$ such that $x_i\cdot t\not\in
\cup^n_{i=1} p_i(\frac{r}{4})$, where
$p_i(r):=A_i(V_i,r)$ is the flow neighborhood of $p_i$. \\
\indent There are two cases. The first case is that after meeting
finite flow neighborhoods $p_i(r)$, $x_i\cdot t$
becomes a flow line in $\G_2\cup\G_3$. Then \\
$$
\int_{[x_i,x_i\cdot t]}\a\ge [\frac{t-\tau}{T_0}]\r-C,
$$
where $\tau,C>0$ are constants. Hence we can take the trajectory
$[x_i,x_i\cdot t]$ as the required $(\mathcal{U}(s_i), T)$-chain
for $s_i=0$ and $T$ arbitrary large. The second case is that
$x_i\cdot t$ can meet infinitely many flow neighborhoods
$\{p_{i_0}(r),\cdots,p_{i_j}(r),\cdots\}$, then
$$
\int_{[x_i,x_i\cdot t]}\a\ge
\sum^{l(t)-1}_{k=0}\int_{[x_i,x_i\cdot
t]_k}\a-\sum^{l(t)-1}_{k=1}\osc_{x\in p_{l_k}(r)}\b_{p_{l_k}}(x),
$$
where $l(t)$ is an integer representing the number of flow
neighborhoods that $[x_i,x_i\cdot t]$ has met and $[x_i,x_i\cdot
t]_k\in \G_1$ is a trajectory of $[x_i,x_i\cdot t]$ between
$p_{i_k}(r)$ and $p_{i_{k+1}}(r)$. We have
$$
\int_{[x_i,x_i\cdot t]}\a\ge l(t)(1-\la)\r.
$$
If $t$ is large enough, this integral is larger than any
given number. Hence we can take $[x_i,x_i\cdot t]$ as the required chain.\\

(2). For any $i$, $x_i\cdot t\in \cup^n_{i=1}B_{\frac{r}{4}}(p_i)$
for some $t$.\\
\indent Firstly there exists a $T(\frac{r}{4})$ such that for any
trajectory $\g\in \G_1(\frac{r}{4})$ with time interval $[a,b]$,
where $T(\frac{r}{4})>b-a$ (one can use the same argument in the
proof of Theorem 3.2). Consider the homeomorphism
$\phi_{T(\frac{r}{4})}:X\rightarrow X$ which maps $x$ to $x\cdot
T(\frac{r}{4})$. Let $B_s$ be any balls in $\mathcal{U}(s)$, then
by compactness the upper bound of the diameter of
$\phi_{t}(B_s),\forall t\in [0,T(\frac{r}{4})]$, will be less than
a uniform constant $C_s<\frac{r}{4}$ if $s$ is small enough.

Since $x_0\cdot t_0$ will meet the first flow neighborhood
$p_{i_0}(\frac{r}{4})$, we can replace $x_0$ by $x_0\cdot t_0$,
since $x_0\cdot t_0$ is also a point in the chain recurrent set of
the flow $v$. So we assume that $x_0\in p_{i_0}(\frac{r}{4})$. Now
we want to modify the $(\mathcal{U}(s), T)$-chain (assume
$T>T(\frac{r}{4})$) ${\g}_i=\{x_0,\cdots, x_k=x_0|t_0,\cdots,
t_{k-1}\}$ to another $(\mathcal{U}(r), T(\frac{r}{4}))$-chain
$\g'=\{x'_0,\cdots, x'_k=x'_0|t'_0,\cdots, t'_{k-1}\}$ satisfying
that the pair $(x'_j\cdot t'_j,x'_{j+1})\in p_{i_j+1}(r)$ for
$j=0,1,\cdots,k-1$. If $B_s(m_{1})$ is contained in some
$p_{i_1}(r)$, then we are done. Assume that $B_s(m_{1})$ is not
contained in any $p_{i}(r)$. By the choice of $T(\frac{r}{4})$,
the flow line $[x_1,x_1\cdot T(\frac{r}{4})]$ will intersect some
flow neighborhood $p_{i_1}(\frac{r}{4})$ at some point
$x'_1=x_1\cdot \tau_1\in B_{\frac{r}{4}}(p_{i_1})$. If $s$ is
small enough, then $(x_0\cdot (t_0+\tau_1)\in
B_{\frac{r}{2}}(p_{i_1})$, hence in $p_{i_1}(r)$. Let $x'_0=x_0$,
$x'_1=x_1\cdot \tau_1$, $t'_0=t_0+\tau_1$ and $t'_1=t_1-\tau_1$.
Hence we get the modified chain $\{x'_0, x'_1, x_2,\cdot,
x_k|t'_0,t'_1,t_2,\cdots\}$. Now $x'_1\in p_{i_1}(\frac{r}{4})$,
we can continue this operation until we get the modified chain
$\g'$. If $s$ is small, it is easy to see
\begin{equation*}
\int_\g\a=\int_{\g'}\a\ge \sum^{k-1}_{j=0}(\int_{[x'_j,x'_j\cdot
t'_j]}-\osc_{x\in p_{i_j}(r)}\b_{p_{i_j}}(x))\ge k(1-\la)\r>0
\end{equation*}
So if $n$ is large enough, the integral of $\a$ along $n\g$ will
be larger than any
given number. So (\ref{3.9}) holds.\qed\\
 \\
\indent In view of Theorem 3.2-3.4, the following corollary is
obvious.\\
 \\
{\em {\bf Corollary 3.5}$\;\;$ Let $v$ be an $\a$-flow on the
compact metric space $(X,d)$. If for any chain $\tilde{\g}$,
$$
\left|\int_{\tilde{\g}}\a\right|\le M
$$
for some $M>0$, then $v$ is a gradient-like flow.}\\
 \\
{\bf Example 3.2}\quad Let $v$ be a flow on $S^1$ which has three
fixed points at $\ta=0,\frac{2\pi}{3},\frac{4\pi}{3}$, and the
forward direction of the flow is the anticlockwise direction. It
is easy to see that $v$ carries a nontrivial cocycle and is a non
gradient-like flow.\\
 \\
{\bf Example 3.3}\quad Let $v$ be the flow generated by a closed
Morse 1-form $\o$ on a manifold $X$. Assume $[\o]$ is a nontrivial
cohomology class, then there are some cases:
\begin{description}
\item[(1)] $\o$ has no zero point. So any trajectory is the one we
want to find in Theorem 3.4, hence $v$ is a non gradient-like
flow.

\item[(2)] $\o$ has only one zero point $p$. Then either
$\G_1(r)\neq\emptyset$ or $\G_2(r)\cup \G_3(r)\neq\emptyset$. In
the first case, there is a homoclinic orbit. In the second case,
we can prove by the same argument as in Example 2.1 that for
$\g\in\G_2(r)\cup \G_3(r)$, the integral $\int_{\g(T)}\o$ will be
larger than any given number. Thus by Theorem 3.4, $v$ is a non
gradient-like flow.

\item[(3)] $\o$ has several zero points with the same index. If
$\text{dim}X=1$, and the orbits between them forms a heteroclinic
chain, then $v$ is non gradient-like. Otherwise, it is a
gradient-like flow as shown in Example 3.1. If $\text{dim}X\ge 2$,
then there are also some cases: (i) if there exists a homoclinic
orbit; (ii) if there exists a heteroclinic cycle; (iii) if
$\G_2(r)\cup\G_3(r)\neq\emptyset$. All these three cases imply $v$
is a non gradient-like flow.

\item[(4)] $\text{dim}X \ge 2$, and there are at least two zero
points with different indices. Assume $p,q$ are the two among
them. Let $X^u_p$ be the unstable manifold at $p$ and $X^s_p$ be
the stable manifold at $p$. Denote $p^+=\text{dim}X^u_p$ and
$p^-=\text{dim}X^s_p$. We have the relation
$$
p^++p^-=\text{dim}X,\;q^++q^-=\text{dim}X.
$$
Then either $p^++q^->\text{dim}X$ or $p^-+q^+>\text{dim}X$. In
either cases, there exists trajectories $\g_i$ connecting the two
zero points such that $\int_{\g_i}\o\rightarrow \infty$ as $i\to
\infty$. Hence $v$ is a non gradient-like flow.

\end{description}
So in summary, if $\text{dim}X\ge 2$, the flow $v$ generated by a
nontrivial closed Morse 1-form $\o$ is a non gradient-like flow.\\
\\
{\bf 4. An algebraic theorem relative to $\pi$-Morse decomposition}\\
 \\
\indent In this section, we let $X$ be an $m$-dimensional compact
polyhedron with metric.
Let $v$ be a $\a$-flow on $X$, where $\a$ is a nontrivial cocycle with higher rank.\\
\indent By Theorem 3.4.4 in \cite{FJ2}, the flow $v$ has a
$\pi$-Morse decomposition, where $\pi$ is the deck transformation
group on the covering space $\bar{X}$ determined by the cocycle
$\a$. Assume that $\pi$ is spanned by $G_+:=\{l_1,\cdots,l_s\}$,
where $s$ is the rank of $\a$. By the proof of Theorem 3.4.3 in
\cite{FJ2}, we can get $s$ pair cross sections $N_i^\pm$ that
satisfy the relations:
$$
l_i\cdot N_i^-=N_i^+, i=1,\cdots, s.
$$
The union $\cup_{i=1}^s (N_i^+\cup N_i^-)$ forms the boundary of a
fundamental domain $X_0$. Since $X$ is a compact polyhedron, we
can take an $m-1$-dimensional polyhedron approximating the cross
sections $N_i^-$ for $i=1,\cdots, s$ such that the intersections
$N_\a^-=N_{i_1}^-\cap N_{i_2}^-\cap\cdots N_{i_r}^-,
\a=\{i_1\cdots,i_r\}$ are polyhedra with codimension $r$. So
without loss of generality, we can assume the cross sections
$N_i^\pm$ and their intersections are polyhedra. For convenience,
we denote $N_i^-=N_i$.\\

\
\begin{figure}
\begin{center}
\includegraphics{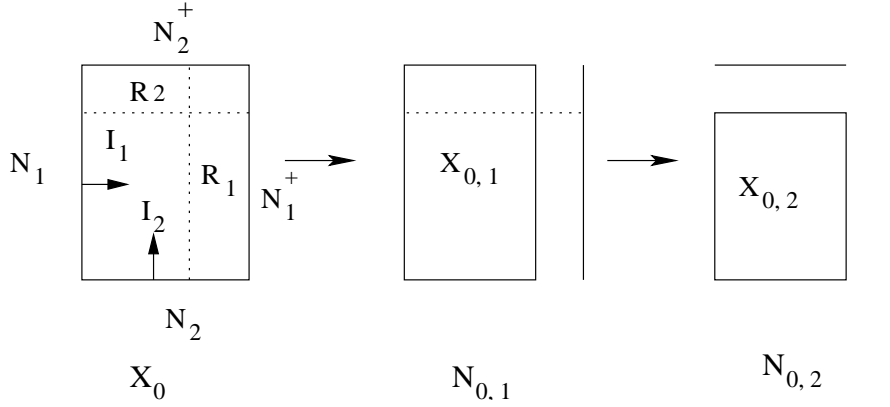}
\caption{Cell decomposition}
\end{center}
\end{figure}

\noindent {\bf Cell decomposition of the fundamental domain}\quad
We can retract $N^+_i$ a small distance $\epsilon$ into the inner
part of the fundamental domain to get a homeomorphic codimension
$1$ polyhedron (see Figure 1. for $s=2$ case). We denote the
obtained polyhedron as $R_i$. Let $i_{j,+}: R_j\rightarrow
X_0,\,j=1,\cdots,s$ be the inclusion and $i_{j,-}: R_j\rightarrow
N_j$ be the homeomorphism.

Let $I=(0,1)$ be the 1-dimensional cell. Firstly we delete
$R_1\times \epsilon I^1$ and get a polyhedron $N_{0,1}$. Denote
$X_{0,1}=N_{0,1}/N_1^+$. We call $G_1=(X_{0,1}, N_{0,1}, R_1,
i_{1,\pm})$ to be the first gluing data. Assume after $r$-th
excision, one has the gluing data $G_r=(X_{0,r}, N_{0,r},
R_1,\cdots, R_r, i_{1,\pm},\cdots, i_{r,\pm})$. Now deleting
$R_{r+1}\times \epsilon I^{r+1}$ from $N_{0,r}$, we get the
polyhedron $N_{0,r+1}$. Define
$X_{0,r+1}=N_{0,r+1}-\cup_{i=1}^{r+1} N_i^+$. Hence after $r+1$-th
excision, one has the gluing data $G_{r+1}$. Continue this process
until we get the $s$-th gluing data $G_s=(X_{0,s}, N_{0,s},
R_1,\cdots, R_s, i_{1,\pm},\cdots, i_{s,\pm})$. $G_i$ is called
the $i$-th gluing data because starting from these data, one can
reconstruct and give a compatible cell decomposition to the
fundamental domain $X_0$ through $i$ times inverse attaching
operations.

Now we want to construct a compatible regular cell decomposition
of $X_0$. Firstly we chose an arbitrary cell decomposition of
$\cup^s_{i=1}(R_i^\circ\cap X_{0,s})$, where $R_i^\circ$ is the
interior part of $R_i$. Then the homeomorphisms
$i_{j,-},j=1,\cdots,s$ provide a cell decomposition of the other
half part of the boundary of $X_{0,s}$. In the interior of
$X_{0,s}$, we choose an arbitrary compatible regular cell
decomposition. Then we get a cell decomposition of $X_{0,s}$. Now
the action of the group provides a cell decomposition of
$N_{0,s}$. At this moment, we can use our gluing data to add the
cylindrical cells to $N_{0,s}$. We add firstly the cell $e\times
\epsilon I^s$ for $e\in R_s\cap N_{0,s}$ to $N_{0,s}$ by the
inclusion map $i_{s,+}$ and by identifying the end $e\times
\epsilon \{1\}$ to $l_s\cdot i_{s,-}(e)$. So we obtain a cell
decomposition of $N_{0,s-1}$ with the gluing data $G_{s-1}$.
Continuing this process, we can obtain a compatible regular cell
decomposition (i.e., each closed sets appeared here are all cell
subcomplexes) of $X_{0}$.

Since $N_{0,s}$ is a cell complex, we can choose arbitrary
orientation of a cell $e\in C_q(N_{0,s})$. Hence fixing the
incidence coefficients appearing in the boundary operator of the
cell chain complex. For the orientation of the cylindrical cells
$e\times I^{i_1}\times \cdots I^{i_k}$, we choose the natural
orientation, i.e., the order $i_1<i_2\cdots<i_k$ defines the
positive orientation. The orientation of the product cell $e\times
I^{i_1}\times \cdots I^{i_k}$ is uniquely determined by the
orientation of the factor cells. So the incidence coefficients of
the boundary operator of $C_*(X_0)$ are uniquely determined.\\
 \\
\noindent \textbf{Deformation of cell chain complexes}\quad Let
$\a\in C_q(X_0)$, then there is a unique decomposition
$$
\a=(-1)^{q-1}(i_{1,+})_*(\a_1)\times \epsilon I^1+\a'_1.
$$
Define $\b_1(\a)=\a_1$, $\b'_1(\a)=\a'_1$ and
$\ta_1(\a_1)=(-1)^{q-1}(i_{1,+})_*(\a_1)\times \epsilon I^1$. We
have the relation
$$
(\begin{array}{cc} \ta_1,&1\end{array})\left(\begin{array}{c}
\b_1\\\b'_1\end{array}\right)=Id.
$$
Similarly, for $\a\in C_q(N_{0,k-1})$ we have the decomposition
$$
\a=(-1)^{q-1}(i_{k,+})_*(\a_k)\times \epsilon I^k+\a'_k.
$$
Let $\b_k(\a):=\a_k\in C_{q-1}(R_k-\cup^{k-1}_{i=1}R_i\times I^i)$
, $\b'_k(\a):=\a'_k\in C_{q}(N_{0,k})$ and
$\ta_k(\a_k):=(-1)^{q-1}(i_{k,+})_*(\a_k)\times \epsilon I^k$. We
have the equalities
\begin{equation}
(\begin{array}{cc} \ta_k,&1\end{array})\left(\begin{array}{c}
\b_k\\\b'_k\end{array}\right)=Id, i=2,\cdots, s\tag{4.1}
\end{equation}\\
 \\
{\em {\bf Lemma 4.1}\quad $\ta_k,\b_k,\b'_k$ are cellular chain
maps and satisfy the following equalities:
\begin{enumerate}
\item[(1)] $d_{N_{0,k-1}}\circ \ta_k+\ta_k\circ
d_{R_k}=l_k\cdot(i_{k,-})_*-(i_{k,+})_*$

\item[(2)] $\b_k\circ d_{N_{0,k-1}}+d_{R_k}\circ \b_k=0$

\item[(3)] $\b'_k\circ d_{N_{0,k-1}}-d_{N_{0,k}}\circ
\b'_k=(l_k\cdot(i_{k,-})_*-(i_{k,+})_*)\circ\b_k$

\end{enumerate}
The above relations can also be written in matrix form:
\begin{enumerate}
\item[(4)] $$ d_{N_{0,k-1}}\circ
(\ta_k,1)=(\ta_k,1)\left(\begin{array}{cc} -d_{R_k}&0\\
f_k& d_{N_{0,k}}\end{array}\right)
$$
\item[(5)] $$ \left(\begin{array}{c}\b_k\\ \b'_k\end{array}\right)
\circ
d_{N_{0,k-1}}=\left(\begin{array}{cc} -d_{R_k}&0\\
f_k& d_{N_{0,k}}\end{array}\right)\left(\begin{array}{c}\b_k\\
\b'_k\end{array}\right)
$$
Here $f_k:=l_k\cdot(i_{k,-})_*-(i_{k,+})_*$.
\end{enumerate}
 }\
\noindent Proof
\begin{enumerate}
\item[(1)] Let $\a\in C_{q-1}(R_k-\cup^{k-1}_{i=1}R_i\times I^i)$,
then
\begin{align*}
&d_{N_{0,k-1}}\circ
\ta_k(\a)=d_{N_{0,k-1}}((-1)^{q-1}(i_{k,+})_*(\a)\times \epsilon
I^k)\\
&=(-1)^{q-1}(i_{k,+})_*(d_{R_k}\a)\times \epsilon
I^k+(l_k\cdot(i_{k,-})_*-(i_{k,+})_*)(\a)\\
&=-\ta_k\circ d_{R_k}\a+(l_k\cdot(i_{k,-})_*-(i_{k,+})_*)(\a)
\end{align*}
\item[(2)] Assume $\a\in C_{q-1}(N_{0,k-1})$ has the form
$$
\a=(-)^{q-2}(i_{k,+})_*(e_1)\times \epsilon I^k+e_2,
$$
where $e_1\in C_{q-2}(R_k-\cup^{k-1}_{i=1}R_i\times \epsilon I^i)$
and $e_2\in C_{q-1}(N_{0,k})$. Then
$$
d_{N_{0,k-1}}\a=(-)^{q-2}(i_{k,+})_*(d_{R_k}e_1)\times \epsilon
I^k+(l_k\cdot(i_{k,-})_*-(i_{k,+})_*)(e_1)+d_{N_{0,k}}e_2,
$$
hence
$$
\b_k\circ d_{N_{0,k-1}}\a=-d_{R_k}\circ \b_k(\a),
$$
and
$$
\b'_k\circ d_{N_{0,k-1}}\a=(l_k\cdot(i_{k,-})_*-(i_{k,+})_*)\circ
\b_k(\a)+d_{N_{0,k}}\circ \b'_k(\a)
$$
\end{enumerate}
\qed\\
 \\
{\em {\bf Lemma 4.2 }\quad Let $f_k:=l_k\cdot
(i_{k,-})_*-(i_{k,+})_*,k=1,\cdots,s$ be the chain map from
$C_{*-1} (R_k-\cup^{k-1}_{i=1}R_i\times I^i)$ to
$C_{*-1}(N_{0,k})$. Then $(\ta_k,1)$ is a chain isomorphism from
the cone $con(f_k):=(C_{*-1} (R_k-\cup^{k-1}_{i=1}R_i\times
I^i)\oplus C_{*-1}(N_{0,k}), d_{c}^k)$ to $(C_*(N_{0,k-1}),
d_{N_{0,k-1}})$, where the differential of the algebraic cone is
$$
d_{c}^k:=\left(\begin{array}{cc} -d_{R_k}& 0\\
f_k& d_{N_{0,k}}\end{array}\right)
$$
}\\
 \\
Proof. By equalities (4) and (5) of Lemma 4.1, $(\ta_k,1)$
and $\left(\begin{array}{c}\b_k\\
\b'_k\end{array}\right)$ are chain maps with respect to
corresponding cellular chain complexes. Furthermore, we have
$$
(\ta_k,1)\left(\begin{array}{c}\b_k\\
\b'_k\end{array}\right)=Id
$$
and
$$
\left(\begin{array}{c}\b_k\\
\b'_k\end{array}\right)(\ta_k,1)=\left(\begin{array} {cc} 1&0\\
0&1\end{array}\right).
$$
Therefore $(\ta_k,1)$ is a chain isomorphism.\qed\\
 \\
{\em {\bf Lemma 4.3}\quad Let
$$
D^k_q=C_{q-1}(R_1)\oplus C_{q-1}(R_2-(i_{1,+})(R_1)\times \epsilon
I^1)\oplus \cdots\oplus
C_{q-1}(R_k-\cup^{k-1}_{j=1}(i_{j,+})(R_j)\times \epsilon
I^j)\otimes C_q(N_{0,k})
$$
\begin{equation}\label{lm-4.3}
d_{c,k}=\left(\begin{array}{ccccc} -d_{R_1}&\cdots&0&0&0\\
\vdots&\ddots&\vdots&\vdots&\vdots\\
f^k_{k-1,1}&\cdots&-d_{R_{k-1}}&0&0\\
f^k_{k,1}&\cdots&f^k_{k,k-1}&-d_{R_k}&0\\
{f^k}'_{k,1}&\cdots&{f^k}'_{k,k-1}&{f^k}'_{k,k}&d_{N_{0,k}}
\end{array}\right)\tag{4.2}
\end{equation}
Here
$$
f^k_{m,n}=\b_m\circ \b'_{m-1}\cdots\b'_{n+1}\circ f_{n},
$$
for $1\le n<m\le k$;
$$ {f^k}'_{k,n}=\b'_k\circ
\b'_{k-1}\cdots \b'_{n+1}\circ f_{n},
$$
for $1\le n\le k-1$; and
$$
{f^k}'_{k,k}=f_k.
$$
Then $(D^k_*, d_{c,k})$ is a chain complex. }\\
 \\
Proof. We prove $d_{c,k}\circ d_{c,k}=0$ by using induction for
$k$. If $k=1$, then $D^1_*=C_{*-1}(R_1)\oplus C_*(N_{0,1})$ and
$$
d_{c,1}=\left(\begin{array}{cc}-d_{R_1}&0\\
f_1& d_{N_{0,1}}\end{array}\right).
$$
It is obviously a chain complex. Assume that $(D^k_*,d_{c,k})$ is
a chain complex, we want to show that $(D^{k+1}_*, d_{c,k+1})$ is
a chain complex. Notice that for $1\le n<m\le k$, we have
$$
f^{k+1}_{m,n}=f^k_{m,n}.
$$
Therefore $d_{c,k}$ and $d_{c,k+1}$ can be expressed by block
matrices:
$$
d_{c,k}=\left(\begin{array}{cc}\tilde{d}_k&0\\
\tilde{f}_k &d_{N_{0,k}}\end{array}\right)
$$
and
$$
d_{c,k+1}=\left(\begin{array}{ccc}\tilde{d}_k&0&0\\
\b_{k+1}\cdot\tilde{f}_k
&-d_{R_{k+1}}&0\\
\b'_{k+1}\cdot\tilde{f}_k&f_{k+1}&d_{N_{0,k+1}}\end{array}\right)
$$
According to the assumption, we have the relations:
$$
\tilde{d}_k\circ \tilde{d}_k=0,
$$
and
$$
\tilde{f}_k\tilde{d}_k+d_{N_{0,k}}\tilde{f}_k=0.
$$
Let $\b_{k+1}$ and $\b'_{k+1}$ act on the above equality and using
the formulas (2) and (3) in Lemma 4.1, we obtain
$$
\b_{k+1}\cdot \tilde{f}_k\cdot \tilde{d}_k-d_{R_{k+1}}\cdot
\b_{k+1}\cdot\tilde{f}_k=0,
$$
and
$$
\b'_{k+1}\cdot \tilde{f}_k\cdot \tilde{d}_k+d_{N_{0,k+1}}\cdot
\b'_{k+1}\cdot\tilde{f}_k+f_{k+1}\cdot\b_{k+1}\tilde{f}_k=0.
$$
combining the above two relations, the fact that $\tilde{d}_k\circ
\tilde{d}_k=0$ and Lemma 4.2, we can deduce $d_{c,k+1}\cdot d_{c,k+1}=0$.\qed\\
 \\
{\em {\bf Lemma 4.4}\quad $\forall k=1,\cdots,s$, the cellular
chain complex $(D^k_*,d_{c,k})$ is chain isomorphic to
$(C_*(X_0),d_{X_0})$.}\\
 \\
Proof. Firstly we will show the following map
$$
\left(\begin{array}{ccccc} 1&\cdots&0&0&0\\
\vdots&\ddots&\vdots&\vdots&\vdots\\
0&\cdots&1&0&0\\
0&\cdots&0&\ta_{k+1}&1\end{array}\right)_{(k+1)\times(k+2)}:C_{q-1}(R_1)\oplus
C_{q-1}(R_2-(i_{1,+})(R_1)\times \epsilon I^1)\oplus \cdots
$$
\begin{align*}
&\oplus C_{q-1}(R_{k+1}-\cup^{k}_{j=1}(i_{j,+})\times \epsilon
I^j)\otimes
C_q(N_{0,k+1})\rightarrow \\
&C_{q-1}(R_1)\oplus C_{q-1}(R_2-(i_{1,+})(R_1)\times \epsilon
I^1)\oplus \cdots\oplus
C_{q-1}(R_{k}-\cup^{k-1}_{j=1}(i_{j,+})\times \epsilon I^j)\otimes
C_q(N_{0,k})
\end{align*}
 is a chain map, i.e., it should satisfy
the commutation relation:
$$
\left(\begin{array}{ccccc} 1&\cdots&0&0&0\\
\vdots&\ddots&\vdots&\vdots&\vdots\\
0&\cdots&1&0&0\\
0&\cdots&0&\ta_{k+1}&1\end{array}\right)\circ
d_{c,k+1}=d_{c,k}\circ\left(\begin{array}{ccccc} 1&\cdots&0&0&0\\
\vdots&\ddots&\vdots&\vdots&\vdots\\
0&\cdots&1&0&0\\
0&\cdots&0&\ta_{k+1}&1\end{array}\right).
$$
This is equivalent to the following equality:
$$
\left(\begin{array}{ccc}\tilde{d}_k&0&0\\
\ta_{k+1}\b_{k+1}\tilde{f}_k+\b'_{k+1}\cdot\tilde{f}_k&-\ta_{k+1}d_{R_{k+1}}+f_{k+1}&d_{N_{0,k+1}}\end{array}\right)=
\left(\begin{array}{ccc}\tilde{d}_k&0&0\\
\tilde{f}_k&d_{N_{0,k}}\ta_{k+1}&d_{N_{0,k}}\end{array}\right),
$$
where $\tilde{f}_k,\tilde{d}_k$ are defined in Lemma 4.3. Now this
is true by Lemma 4.1 and the fact
$d_{N_{0,k}}|_{N_{0,k+1}}=d_{N_{0,k+1}}$. Thus we proved
$$
\left(\begin{array}{ccc} Id&0&0\\
0&\ta_{k+1}&1\end{array}\right)
$$
is a chain map. It is also an isomorphism, since it has an inverse
chain map:
$$
\left(\begin{array}{cc} Id&0\\
0&\b_{k+1}\\
0&\b'_{k+1}\end{array}\right).
$$
Therefore we proved $(D^{k+1}_*,d_{c,k+1})$ is chain isomorphic to
$(D^k_*,d_{c,k})$ for any $k=1,\cdots,s$. Note the $k=1$ case was
already proved in Proposition 4.1.1 in \cite{FJ2}.
\qed\\
 \\
Define $\pi_+$ as the monoid constructed by $(G_+,e)$ with the
group action from $\pi$. Let $\r_e:\pi_1(X)\to \pi$ be the
extension of $\pi$ by the normal group $\pi_1(\bar{X})$. Since
$\pi_+$ is a monoid in $\pi$, the set $\r^{-1}_e(\pi_+)$ is also a
monoid of $\pi_1(X)$. We denote $\r^{-1}_e(\pi_+)$ by
$\pi_1(X)_+$. ${\Bbb Z}\pi_1(X)_+$ is a subring of ${\Bbb
Z}\pi_1(X)$.

Tensor $(D^k_*,d_{c,k})$ with the ring ${\Bbb Z}\pi_1(X)_+$, then
we have the ${\Bbb Z}\pi_1(X)_+$-module chain complex $({\Bbb
Z}\pi_1(X)_+\otimes D^k_*,I\otimes d_{c,k})$.

On the other hand, we can lift the related quantities
(subcomplexes, or maps) to the universal covering space
$\tilde{X}$ such that they become $\pi_1(X)_+$-equivariant. For
instance, we can lift the $m-1$-dimensional polyhedra $\cup_{g\in
\pi_+} g\cdot R_i$ to $\tilde{X}$ such that the obtained polyhedra
$\tilde{R}_{i}$ is $\pi_1(X)_+$-equivariant. If $F$ is a map
($\b_k,\b'_k,f_k$, etc.), then the action of the lifted map
$\tilde{F}$ is defined as: $\tilde{F}(g\cdot e)=g\cdot
\tilde{F}(e)$. Thus we can get a ${\Bbb Z}\pi_1(X)_+$-module
cellular chain complexes $(\tilde{D}^k_*,\tilde{d}_{c,k})$, where
$$
\tilde{D}^k_q:=C_{q-1}(\tilde{R}_1)\oplus
C_{q-1}(\tilde{R}_2-(\tilde{i}_{1,+})(\tilde{R}_1)\times \epsilon
I^1)\oplus \cdots\oplus
C_{q-1}(\tilde{R}_k-\cup^{k-1}_{j=1}(\tilde{i}_{j,+})\times
\epsilon I^j)\oplus C_q(\tilde{N}_{0,k})
$$
and $\tilde{d}_{c,k}$ is defined by the same form of the matrix
(\ref{lm-4.3}), and the difference is that the entries of
(\ref{lm-4.3}) are replaced by their lifting homomorphisms.

It is easy to see that the two ${\Bbb Z}\pi_1(X)_+$-module chain
complexes $({\Bbb Z}\pi_1(X)_+\otimes D^k_*, I\otimes d_{c,k})$
and $(\tilde{D}^k_*,\tilde{d}_{c,k})$ are chain equivalent.\\

Let $\tilde{X}_+=\cup_{g\in \pi_1(X)_+} g\cdot
\tilde{X}_0$. By Lemma 4.4, it is obvious that the following holds.\\
 \\
{\em {\bf Proposition 4.5} The two ${\Bbb Z}\pi_1(X)_+$-module
chain complexes $(C_*(\tilde{X}_+),d_{\tilde{X}_+})$ and
$(\tilde{D}^s_*,\tilde{d}_{c,k})$ are chain equivalent. }\\
 \\
\textbf{Monodromy representations and Novikov numbers}\quad In
this part, we will generalize the crucial theorem, Theorem 4.2.1
in \cite{FJ2} for an integral cocycle $\a$ to a higher rank
cocycle $\a$. To get such a theorem, one needs also to study the
monodromy
representation and the evaluation representation as in \cite{FJ2}.\\
 \\
Since $\pi$ is a commutative group with rank $s$, we obtain a
group homomorphism
$$
\r_{\pi_1}:\,\pi_1(X)\xrightarrow{\r_e} \pi\xrightarrow{\r_A}
\mathbb{Z}^s.
$$
Here $\r_A(l_1^{a_1}\cdots l_s^{a_s}):=(a_1,\cdots,a_s)$. Thus we
have a ring homomorphism
$$
\r_{\pi_1}:\,\mathbb{Z}[\pi_1(X)]\longrightarrow\mathbb{Z}[\mathbb{Z}^s].
$$
This ring homomorphism can induce another ring homomorphism
$$
\r_q:\,\mathbb{Z}[\pi_1(X)]\rightarrow
Q_s=\mathbb{Z}[t_i,t_i^{-1};\,i=1,2,\cdots,s]
$$
defined as, for $g=\sum z_j (l_1^{a_{j_1}}\cdot l_s^{a_{j_s}})\in
\mathbb{Z}[\pi]$,
$$
\r_q(g)=\sum z_j t_1^{a_{j_1}}\cdots t_s^{a_{j_s}}.
$$
In fact, $\r_{q}$ is fully determined by the group $\pi$ and its
representation. Restricting $\r_{q}$ to the subring
$\mathbb{Z}[\pi_1(X)_+]$, we can get a ring homomorphism
$$
\r_p=\r_{q}|_{\mathbb{Z}[\pi_1(X)_+]}:\,\mathbb{Z}[\pi_1(X)_+]\longrightarrow
P_s=\mathbb{Z}[t_1,\cdots,t_s]
$$\\
 \\
Let $\tilde{E}$ be a local system of free abelian groups on the
compact polyhedron $X$, then $\tilde{E}$ is determined by its
monodromy representation $\r_{\tilde{E}}$:
$$
\r_{\tilde{E}}:\,\pi_1(X,x_0)\longrightarrow
\text{Aut}(\tilde{E}_0)=GL(k,\mathbb{Z})
$$
where $\tilde{E}_0$ is the fibre of the free abelian group at
$x_0$ and $k=\text{rank}(\tilde{E}_0)$. Let
$E=\tilde{E}\otimes\mathbb{C}$, then $E$ is a complex flat vector
bundle with the holonomy $\r_E$
$$
\r_E:\,\pi_1(X,x_0)\longrightarrow
GL(k;\mathbb{Z})\otimes\mathbb{C}
$$
\indent
Now the tensor product of the representations $\r_q\otimes\r_{\tilde{E}}$
gives a representation of a $\mathbb{Z}[\pi_1(X)]$-ring to the linear space
$(Q_s)^k$, where $Q_s$ is the polynomial space with $s$ variables over $\mathbb{Z}$.\\
\indent
Since $\r_{\tilde{E}}$ is an anti-homomorphism, i.e., $\forall g,g'\in \mathbb{Z}[\pi_1(X)],
 \r_E(g\cdot g')=\r_E(g')\cdot\r_E(g)$, hence $\r_P\otimes\r_{\tilde{E}}$ gives a
 right $\mathbb{Z}[\pi_1(X)]_+$-module structure on $P^k_s$. With the $P_s$-module
 structure of itself, $P_s^k$ becomes a $(P_s,\mathbb{Z}[\pi_1(X)_+])$-bimodule.\\
\indent
Define $D_*=P_s^k\otimes_{\mathbb{Z}[\pi_1(X)]_+}C_*(\tilde{X}_+)$, then $D_*$ is a $P_s$-module chain complex.\\
 \\
\textbf{ Evaluation representations}\quad Take any complex
$s$-vector $a=(a_1,\cdots, a_s)\in \mathbb{C}^s$. The complex
number field $\mathbb{C}$ can be given a $P_s$-module structure,
whose module structure is provided by the action: for a polynomial
$P(t_1,\cdots,t_s)$, $P(t_1,\cdots,t_s)\cdot
x=P(a_1,\cdots,a_s)\cdot x=P(a)\cdot x$ for $x\in \mathbb{C}$. We
denote the $P_s$-module of $\mathbb{C}$ evaluated at $t=a$ by
$\mathbb{C}_a$. Similarly for any $a\in
(\mathbb{C}^*)^s,\,\mathbb{C}$ can be viewed as a $Q_s$-module. If
$p$ is a prime number, then the field $\mathbb{Z}_p$ also has a
$P_s$-module structure which is given by the evaluation at $t=0$.
We consider the complexes $\mathbb{C}_a\otimes_{P_s}D_*$ and
$\mathbb{Z}_p\otimes_{P_s}D_*$. The following theorem is the
generalization of Theorem 4.2.1 in \cite{FJ2} for higher rank $\a$.\\
  \\
{\em {\bf Theorem 4.6}\quad Assume $\text{rank}\pi=s\ge 1$.
Let $D_*=P^k_s\otimes_{\mathbb{Z}[\pi_1(X)_+]}C_*(\tilde{X}_+)$ be defined by the above argument. We have \\
\begin{enumerate}
\item[(1)] For any nonzero complex vector $a\in(\mathbb{C}^*)^s$,
the homology $H_*(\mathbb{C}_a\otimes_{P_s}D_*)$ is isomorphic to
$H_*(X,\,a^\pi\otimes E)$, which is viewed as the homology of the
presheaf $a^\pi\otimes E$ on $X$.

\item[(2)] Let $p$ be a prime number and let $\mathbb{Z}_p$ have
the $P_s$-module structure which is provided by the evaluation at
$t=(t_1,\cdots,t_s)=0$. Then the homology
$H_*(\mathbb{Z}_p\otimes_{P_s}D_*)$ is isomorphic to $H_*(X_{0,s},
\cup^s_{j=1}(i_{j,+})(R_j);\,\mathbb{Z}_p\otimes
P_\pi^*\tilde{E})$, where $\tilde{E}$ is a local system on $X$ and
$P_\pi^* E$ is the pull-back local system on $X_{0,s}$ by the
projection $P_\pi:\,\bar{X}\longrightarrow X$.\\

\item[(3)] $H_*(C_0\otimes_{P_s}D_*)$ is isomorphic to
$H_*(X_{0,s},\cup^s_{j=1}(i_{j,+})(R_j);\,P_\pi^*E)$.
\end{enumerate}
}
Proof.\quad Since the proof of (3) is the same as that of (2), we only give the proofs of (1) and (2).\\
\\
(1). Since the $\mathbb{Z}[\pi_1(X)_+]$-basis of
$C_*(\tilde{X}_+)$ is finite, all the complexes related to
$C_*(\tilde{X}_+)$ are finitely generated,
 and hence all the homology groups are finitely generated.\\
\indent For $a\in (C^*)^s$, we have the isomorphism
$$
\aligned
\mathbb{C}_a\otimes_{P_s}D^*\cong& \mathbb{C}_a\otimes_{P_s}((P_s)^k\otimes_
{\mathbb{Z}[\pi_1(X)_+]}C_*(\tilde{X}_+))\\
\cong&(\mathbb{C}_a\otimes_{Q(s)}(Q_s)^k)\otimes_{\mathbb{Z}[\pi_1(X)]}
(\mathbb{Z}[\pi_1(X)]\otimes_{\mathbb{Z}[\pi_1(X)_+]}C_*(\tilde{X}_+))\\
\cong&(\mathbb{C}_a\otimes_\mathbb{Z}\mathbb{Z}^k)\otimes_{\mathbb{Z}[\pi_1(X)]}C_*(\tilde{X})\\
\cong&\mathbb{C}^k\otimes_{\mathbb{Z}[\pi_1(X)]}C_*(\tilde{X}).
\endaligned
$$
Here the representation of $\mathbb{Z}[\pi_1(X)]$ is given by
$$
g\longrightarrow \r_{P(a)}(g)\otimes\r_E(g).
$$
Hence the homology of
$\mathbb{C}^k\otimes_{\mathbb{Z}[\pi_1(X)]}C_*(\tilde{X})$ is the
same as the homology of the presheaf $a^\pi\otimes E$ on $X$ that
corresponds to the flat vector bundle produced by the above holonomy representation. (1) is proved.\\
 \\
(2)
$$
\mathbb{Z}_p\otimes_{P_s}D^*\cong\mathbb{Z}_p\otimes_{P_s}(P^k_s\otimes_{\mathbb{Z}[\pi_1(X)_+]}C_*(\tilde{X}_+))\cong
\mathbb{Z}^k_p\otimes_{\mathbb{Z}[\pi_1(X)_+]}C_*(\tilde{X}_+).
$$
Here the representation of $\mathbb{Z}[\pi_1(X)_+]$ on
$\mathbb{Z}^k_p$ is
\begin{equation}\label{rep-vanish}
g\longrightarrow \r_{P(0)}(g)\otimes\r_E(g).\tag{4.3}
\end{equation}
Since $\r_{P(0)}(g)=\r_{P(0)}\cdot\r_A\r_e(g)$, except in the case
that $g\in\mathbb{Z}[\pi_1(X)_+]$ satisfies $\r_e(g)=1$, the
evaluation representation will make the final representation
vanish. Hence (\ref{rep-vanish}) becomes
$$
\aligned
g\longrightarrow \r_E(g),\;\;&\hbox{if}\;\r_e(g)=1\\
g\longrightarrow 0,\;\;&\hbox{if}\;\r_e(g)\neq 1.
\endaligned
$$
By Proposition 4.5, in order to prove (2), we need to prove that
$\mathbb{Z}^k_p\otimes_{\mathbb{Z}[\pi_1(X)_+]}\tilde{D}^s_*$ is
equivalent to $(\mathbb{Z}_p\otimes P_\pi^*\tilde{E})\otimes
C_*(X_{0,s}, \cup^s_{j=1}(i_{j,+})(R_j))$. Now we have

$$
\mathbb{Z}^k_p\otimes_{\mathbb{Z}[\pi_1(X)_+]}\tilde{D}^s_q=
(\mathbb{Z}^k_p\otimes_{\mathbb{Z}[\pi_1(X)_+]}C_{q-1}(\tilde{R}_1))\oplus\cdots
(\mathbb{Z}^k_p\otimes_{\mathbb{Z}[\pi_1(X)_+]}\oplus
C_q(\tilde{N}_{0,s})).
$$

If we take $\tilde{E}\otimes
C_*(R_t-{\cup^{t-1}_{j=1}(i_{j,+})(R_j)\times \epsilon I^j},N^+)$
as the complexes with twisted coefficients $\tilde{E}$, then
$\mathbb{Z}_p^k\otimes_{\mathbb{Z}[\pi_1(X)_+]}C_*(\tilde{R}_t-
\cup^{t-1}_{j=1}(\tilde{i}_{j,+})(\tilde{R}_j)\times \epsilon
I^j)\cong \mathbb{Z}_p\otimes(\tilde{E}\otimes
C_*(R_t-{\cup^{t-1}_{j=1}(i_{j,+})(R_j)\times \epsilon I^j},N^+))$
(Notice that the evaluation map at $t=0$ makes the cell lying in
$N^+$ vanishing). Similarly
$\mathbb{Z}_p^k\otimes_{\mathbb{Z}[\pi_1(X)_+]}C_*(\tilde{N}_{0,t})\cong
\mathbb{Z}_p\otimes(P_\pi^*\tilde{E}\otimes C_*(X_{0,t}))$.
Therefore
$$
\mathbb{Z}^k_p\otimes_{\mathbb{Z}[\pi_1(X)_+]}\tilde{D}^s_*\cong
(\mathbb{Z}_p\otimes(\tilde{E}\otimes
C_{*-1}(R_1,N^+)))\oplus\cdots\oplus(\mathbb{Z}_p\otimes(P_\pi^*\tilde{E}\otimes
C_*(X_{0,s}))
$$
Denote this complex by $D^{s,0}_*$. Its differential becomes
\begin{equation}\label{bd-op0}
d_{c}^{s,0}=\left(\begin{array}{ccccc} -d_{R_1}&\cdots&0&0&0\\
\vdots&\ddots&\vdots&\vdots&\vdots\\
i^s_{s-1,1}&\cdots&-d_{R_{s-1}}&0&0\\
i^s_{s,1}&\cdots&i^s_{s,s-1}&-d_{R_s}&0\\
{i^s}'_{s,1}&\cdots&{i^s}'_{s,s-1}&{i^s}'_{s,s}&d_{N_{0,s}}
\end{array}\right)\tag{4.4}
\end{equation}
Here
$$
i^s_{m,n}=-\b_m\circ \b'_{m-1}\cdots\b'_{n+1}\circ i_{n,+},
$$
for $1\le n<m\le s$,
$$ {i^s}'_{s,n}=-\b'_s\circ
\b'_{s-1}\cdots \b'_{n+1}\circ i_{n,+},
$$
for $1\le n\le s-1$, and
$$
{i^s}'_{s,s}=-i_{s,+}.
$$

Now for convenience, we assume the coefficient ring is the integer
ring.

Define maps \begin{align*} &I_1:=(i^s_{2,1},\cdots,
i^s_{s,1},i^{s'}_{s,1}):C_{*}(R_1,N^+)\rightarrow \\
&D^{s,1}_*:= C_{*-1}(R_2-(i_{1,+})(R_1)\times \epsilon
I^1,N^+)\oplus\cdots\oplus C_{*}(X_{0,s}), \end{align*} and
$$
P_1:=D^{s,1}_*\rightarrow
\bar{D}^{s,1}_*:=C_{*-1}(R_2-(i_{1,+})(R_1)\times I^1, R_1\cup
N^+)\oplus\cdots\oplus C_{*}(X_{0,s}, R_1)
$$
to be the quotient map. Let $d_c^{s,1}$ be the matrix obtained
from (\ref{bd-op0}) by deleting the first row and the first
column. We have $d_c^{s,1}\circ d_c^{s,1}=0$. Hence
$(D^{s,1}_*,d_c^{s,1})$ is a chain complex. It induces the
quotient chain complex $(\bar{D}^{s,1},\bar{d}_c^{s,1})$. The
following is a short exact sequence.
\begin{equation}\label{sex}
0\rightarrow
C_{*}(R_1)\xrightarrow{I_1}D^{s,1}_*\xrightarrow{P_1}\bar{D}_*^{s,1}\rightarrow
0\tag{4.5}
\end{equation}
To see this, we have to show the following three points:
\begin{enumerate}
\item[(i)] $P_1\circ I_1=0$

\item[(ii)] $\text{ker}P_1\subset \text{im}I_1$

\item[(iii)] $I_1, P_1$ are chain maps.

\end{enumerate}

(i) is obvious. Assume $P_1(\a_2,\cdots,\a_s,\a'_s)=0$, i.e.,
$\forall j=2,\cdots, s$,
$$
\a_j\in R_1\cap R_j-\cup^{j-1}_{t=1} R_t\times \epsilon I^t-N^+
$$
and
$$
\a'_s\in X_{0,s}-\cup^s_{j=2}R_j.
$$
Let
$$
\a=-(\sum^s_{j=2}\ta_j(\a_j)+\a'_s).
$$
Then we have $I_1(\a)=(\a_2,\cdots,\a_s,\a'_s)$. Hence $(ii)$ is
proved. To prove $(iii)$, one need only to check the commutative
diagram:
$$
\begin{CD}
C_q(R_1)@>{I_1}>>D^{s,1}_q\\
@VV{d_{R_1}}V@VV{d^{s,1}_c}V\\
C_{q-1}(R_1)@>{I_1}>>D^{s,1}_{q-1}
\end{CD}
$$
i.e., check the equality:$d^{s,1}_c\cdot I_1=I_1\cdot d_{R_1}$.
This holds in view of the relation $d^{s,0}\cdot d^{s,0}=0$. Hence
$(iii)$ is proved.

From the short exact sequence (\ref{sex}), we obtain the long
exact sequence:
\begin{align}\label{lex1}
&\cdots\rightarrow H_q(R_1)\xrightarrow{I_1}
H_q(D^{s,1}_*)\rightarrow
H_q(\bar{D}^{s,1}_*)\rightarrow\nonumber\\
&\rightarrow
H_{q-1}(R_1)\xrightarrow{I_1}H_{q-1}(D^{s,1}_*)\rightarrow\cdots\tag{4.6}
\end{align}
Similarly, from the short exact sequence
$$ 0 \rightarrow C_*(D^{s,1}_*)
\rightarrow C_*(D^{s,0}_*) \rightarrow (C_{*-1}(R_1),-d_{R_1})
\rightarrow 0
$$
we get the long exact sequence
$$
\cdots \rightarrow H_q(R_1) \xrightarrow{I_1}
H_q(D^{s,1}_*)\rightarrow H_q(D^{s,0}_*)
 \rightarrow H_{q-1}(R_1)\xrightarrow{I_1} H_{q-1}(D^{s,1}_*)
 \rightarrow\cdots
$$
The two long exact sequences yield the following two short exact
sequences
\begin{align*}
0 \rightarrow \text{Coker}(I_1:\,H_q(R_1)\to H_q(D^{s,1}_*)) \rightarrow & H_q(\bar{D}^{s,1}_*)\rightarrow \\
 \rightarrow &\text{Ker}(I_1:\,H_{q-1}(R_1)\to H_{q-1}(D^{s,1}_*)) \rightarrow  0
\end{align*}
and
\begin{align*}
0 \rightarrow\text{Coker}(I_1:\,H_q(R_1)\to H_q(D^{s,1}_*)) \rightarrow & H_q(D^{s,0}_*)\rightarrow \\
\rightarrow &\text{Ker}(I_1:\,H_{q-1}(R_1)\to H_{q-1}(D^{s,1}_*))
\rightarrow 0
\end{align*}
Using the Five-Lemma, we get for any $q\ge 0$
$$
H_q(\bar{D}_*^{s,1})\cong H_q(D^{s,0}_*).
$$
Continuing our operation, finally we get
$$
H_q(D^{s,0}_*)\cong
H_q(\bar{D}^{s,s}_*)=H_q(X_{0,s},\cup^s_{j=1}R_j).
$$
Hence we finally proved the second conclusion of this
theorem.\qed\\
  \\
{\em Novikov numbers}\quad In Theorem 4.4, we have considered the
complex $D_*=P^k_s\otimes_{{\Bbb Z}
 [\pi_1(X)_+]}C_*(\tilde{X}_+)$. In this part, we always let the
vector bundle $E$ that appears in Theorem 4.6 be a trivial line
bundle. Then the complex $D_*$ there has the form
$D_*=P_s\otimes_{{\Bbb Z}[\pi_1(X)_+]}C_*(\tilde{X}_+)$. Since the
representation of ${\Bbb Z}[\pi_1(X)_+]$ in $P_s$ is completely
determined by the cohomology class $[\a]$, we denote the homology
group $H_*({\Bbb C}_a\otimes_{P_s}D_*)$ as $H_*(X,a^\a)$, or in
other words, view the homology group $H_*({\Bbb
C}_a\otimes_{P_s}D_*)$ as the homology group of the presheaf
$a^\a$ on $X$ which is given by the monodromy
representation $\r_P:{\Bbb Z}[\pi_1(X)_+]\to P_s={\Bbb Z}[t_1,\cdots,t_s]$.\\
 \\
\indent
The following is essentially given in \cite{No3}.\\
 \\
{\em {\bf Proposition 4.7}$\;\;$ Let $X$ be a compact polyhedron.
Define a function for fixed $i$ to be
$$
a\in ({\Bbb C}^*)^s\longrightarrow \dim_{{\Bbb C}}H_i(X,a^\a),
$$
then it has the following properties:\\
\indent
(1) It is generically constant, more precisely, except on a
proper algebraic subvariety $\mathcal{L}$ in $({\Bbb C}^*)^s$,
the dimension $\dim_{{\Bbb C}}H_i(X,a^\a)$ is constant and this
constant is just the Novikov number $b_i([\a])$ we defined above.\\
\indent (2) For any point $\hat{a}\in\mathcal{L}$,
$$
\dim_{{\Bbb C}}H_i(X,\hat{a}^\a)>b_i([\a])
$$
}\\
\noindent
{\bf Definition}$\;\; b_i([\a]):=\hbox{rank}(H_i(D_*))$ for $i=0,1,\cdots,m$ are called the Novikov numbers.\\
  \\
{\bf Remark}$\;\;$ If $s=1$ and $X$ is a closed manifold, then it
was proved in \cite{Fa5} that the Novikov number we defined here
is the same as the Novikov number which is defined as
the dimension of the homology group of the Novikov complex. \\
   \\
{\bf 5. Novikov-Morse type inequalities for higher rank $\a$}\\
 \\
\indent Theorem 4.6 is crucial in the proof of Novikov-Morse type
inequalities. As long as we obtain theorem 4.6 for flows carrying
a higher rank cocycle $\a$, we can follow the same line as in the
proof in \cite{FJ2} to get the expected analogous inequalities. We
omit the similar proof as in section 4.4 and 4.5 of \cite{FJ2},
and only
list the results. \\
 \\
{\em Ideals}\quad Let $P_s={\Bbb Z}[t_1,\cdots,t_s]$ be the
polynomial ring with $s$ variables over the integers ${\Bbb Z}$.
Let $\langle t_1,\cdots,t_s\rangle$ be the ideal generated by
$t_1,\cdots,t_s$. Define $I=\{t\in{\Bbb C}^s\,;\,1+\langle
t_1,\cdots,t_s\rangle=0\}$. Let $a\in ({\Bbb C}^*)^s$ be not in
$I$. Define $I_a$ to be the prime ideal in the polynomial ring
$P_s$ consisting of the polynomials vanishing at $a$. By the
choice of $a$, the free terms of all the polynomials $f(t)\in I_a$
are divisible by some prime number $p$. Therefore we obtain
$$
I_a\subset I_p=\langle p\rangle+\langle t_1,\cdots,t_s\rangle
$$
  \\
{\em {\bf Theorem 5.1}$\;\;$ Let $X$ be a compact polyhedron with
a metric $d$. Let $\a$ be an continuous cocycle and $v$ be a
generalized $\a$-flow w.r.t. an isolated invariant set
$\mathcal{A}=\{A_1,\cdots,A_n\}$. Let $\tilde{E}$ be a local
system of free abelian groups and let $E={\Bbb
C}\otimes\tilde{E}$. If $a\in{\Bbb C}^*$ and $a\not\in I=1+\langle
t_1,\cdots,t_s\rangle\subset P_s$, then there is a prime number
$p$ relative to $a$ such that
\begin{align*}
&\sum_{A_i\in \mathcal{A}}p(h(A_i);t,{\Bbb Z}_p\otimes J_{A_i}^* E)=p(X;t,a^\a\otimes E)\notag\\
&+(1+t)Q_1(I_a,I_p;t)+(1+t)Q_2(I_p;t)\tag{$5.1$}
\end{align*}
where $J_{A_i}$ is the inclusion map from the isolating
neighborhood of an isolated invariant set $A_i$ to $X$, and
$Q_1(I_a,I_p;t)$ and $Q_2(I_p;t)$ are all polynomials with
nonnegative integer coefficients.}\\
\\
{\em {\bf Corollary 5.2 (Euler-Poincar\'e formula)}$\;\;$ Under
the hypothesis of Theorem 5.1,
\begin{align*}
\sum_{A\in\mathcal{A}} p(h(A)\,;\,-1\,,\,{\Bbb Z}_p\otimes J^*_A
E)=p(X\,,\,-1\,,\,a^\a\otimes E)\tag{$5.2$}
\end{align*}
In particular, if $E$ is a trivial line bundle, then
\begin{align*}
\sum_{A\in\mathcal{A}}p(h(A)\,;\,-1\,,\,{\Bbb
Z}_p)=\chi(X)\tag{$5.3$}
\end{align*}
for any prime number $p$. Here $\chi(X)$ is the Euler characteristic number of the compact polyhedron.}\\
 \\
\indent As explained in \cite{FJ2}, if the cocycle $\a$ is a
trivial cocycle, then the above Novikov-Morse type inequalities
will induce the Conley-Morse type inequalities given in \cite{CZ}.
If $\o$ is a closed 1-form with higher rank, then Theorem 5.1
generalizes
the rank one case which is proved in \cite{FJ2}. We have\\
 \\
{\em {\bf Theorem 5.3}$\;\;$ Let $\o$ be a closed 1-form with
critical set $\mathcal{A}$ consisting of finitely many connected
components. Let $\tilde{E}$ be a local system of free abelian
groups and $E={\Bbb C}\otimes\tilde{E}$. Assume that $a\in {\Bbb
C}^*$ and $a\not\in I=1+\langle t_1,\cdots,t_s\rangle\subset P_s$.
Then there is a prime $p$ such that
\begin{align*}
&\sum_{A\in\mathcal{A}}p(h(A);t,{\Bbb Z}_p\otimes J^*_A E)\\
&=p(X,t;a^\o\otimes E)+(1+t)Q(E,t)
\end{align*}
where $Q(E,t)$ is a polynomial with nonnegative coefficients.
}\\
 \\
\indent If the closed 1-form $\o$ in Theorem 5.3 has only
non-degenerate isolated zero points or Bott type non-degenerate
zero locus, then the above theorem reduces to the theorems proven
by M. Farber \cite{Fa1,Fa2}(The reader can also
see \cite{FJ2} for a complete description of those cases).\\
 \\
{\em {\bf Corollary 5.4 (Classical Novikov inequality)}$\;\;$ Let
$X$ be an oriented closed smooth manifold and $\o$ a Morse closed
1-form. Then the numbers $c_j(\o)$ of zeros of $\o$ having index
$j$ satisfy
\begin{align*}
c_j(\o)\ge& b_j([\o]) \\
\sum^j_{i=0}(-1)^i c_{j-i}(\o)\ge&\sum^j_{i=0}(-1)^i
b_{j-i}([\o])\tag{$5.4$}
\end{align*}
for $j=0,1,\cdots,m$.\\
}\\
 \\
{\em {\bf Theorem 5.5}$\;\;$ Let $v$ be an $\a$-Morse-Smale flow
w.r.t. $\mathcal{A}=\{A_n,\cdots,A_1\}$. Let $c_j$ be the number
of hyperbolic fixed points with index $j$, $a_j$ be the number of
the hyperbolic periodic orbits with index $j$, and
$\mu_j=c_j+a_j+a_{j+1}$. Then
\begin{align*}
\mu_j\ge& b_j([\a])\\
\sum^j_{i=0}(-1)^{j-i}\mu_i\ge&\sum^j_{i=0}(-1)^{j-i}b_j([\a])\tag{$5.5$}
\end{align*}
for $j=0,1,\cdots,m$.}\\
  \\
{\em {\bf Theorem 5.6 (Vanishing theorem)}$\;\;$ Let $X$ be a
compact polyhedron with a metric $d$. If there exists a flow
carrying a cohomology class $[\a]$ on $X$, then
$$
b_i([\a])=0,\;\;\forall i=0,1,\cdots,m.
$$
}\\
{\bf Example 5.1}\quad Let $f:X\rightarrow X$ be a homeomorphism
from the compact polyhedron $X$ to itself, then we can define the
mapping torus of $f$:
$$
T_f:=\{(x,t)|x\in X, t\in [0.1]\}/(x,0)\sim (f(x),1).
$$
Then the natural projection $\pi:T_f\rightarrow S^1$ provides a
nontrivial cocycle $\a=J(\pi):X\times X\rightarrow {\Bbb R}$
(where $J$ is defined in proposition 3.1.1 of \cite{FJ2}). Define
the flow carrying the cohomology class $\a$ as: for $0\le s<1$,
and $\forall t\in {\Bbb R}$,
$$
[x,s]\cdot t:=[f^{[\![s+t]\!]}(x), (s+t)],
$$
where $[\![\cdot]\!]$ is the function taking the integer part and
$(\cdot)$ is the function taking the decimal part. According to
Theorem 5.6, we have
$$
b_i([\a])=0,\;\forall i.
$$

\vskip 8pt \noindent Huijun Fan, Max-Planck-Institute for
Mathematics in the Sciences, Inselstr. 22-26, 04103 Leipzig, Germany\\
E-mail:\quad hfan@mis.mpg.de\\

\noindent J\"urgen Jost, Max-Planck-Institute for Mathematics in
the
Sciences, Inselstr. 22-26, 04103 Leipzig, Germany\\
E-Mail:\quad jjost@mis.mpg.de\\

\end{document}